\documentclass[11pt, letterpaper, oneside]{myart}
\usepackage{amsmath,amssymb,amscd,amsthm,eucal, graphicx, array}
\usepackage[latin1]{inputenc}
\usepackage{enumitem}

\usepackage[backgroundcolor=white, bordercolor=red, linecolor=red, textsize=tiny]{todonotes} % to do notes
\usepackage{txfonts} % defines \coloneqq (:=)

\usepackage{pifont} % dingbats

\usepackage{color}
\definecolor{linkcol}{rgb}{0, 0, 0.7}
\usepackage[ linkcolor=linkcol, citecolor=linkcol, urlcolor= linkcol, colorlinks=true, hypertexnames=false]{hyperref}

\usepackage{tikz}
\usetikzlibrary{matrix, calc,through,intersections, arrows, decorations.markings}

\renewcommand{\labelenumi}{\theenumi )} % redefines enumerate labels to 1), 2) 3) etc.

\numberwithin{equation}{section}
\renewcommand{\theequation}{\arabic{section}-\arabic{equation}}

\swapnumbers

\theoremstyle{plain}
\newtheorem{theorem}{Theorem}[section]
\newtheorem{lemma}[theorem]{Lemma}
\newtheorem{proposition}[theorem]{Proposition}
\newtheorem{corollary}[theorem]{Corollary}

\theoremstyle{definition}
\newtheorem{definition}[theorem]{Definition}

\newtheorem{notation}[theorem]{Notation}

\newtheorem{note}[theorem]{Note}

\newtheorem{nn}[theorem]{}
     
\newtheorem{igusasarg}[theorem]{Igusa's argument}        
\newtheorem{org}[theorem]{Organization of the paper}  
\newtheorem{twisting}[theorem]{Twisting cochains and twisted tensor products}       
\newtheorem{twistedchain}[theorem]{Twisted chain complex of a fibration}   
\newtheorem{homolfiltr}[theorem]{Homological filtration}   
\newtheorem{mapsfibr}[theorem]{Maps of fibrations}    
\newtheorem{homotfibmap}[theorem]{Maps of homotopy fibers}
\newtheorem{homotfibhomot}[theorem]{Homotopies of maps of homotopy fibers}   
\newtheorem{homotobstr}[theorem]{An obstruction to lifting a homotopy}

\newcommand{\ra}{\rightarrow}
\newcommand{\Ra}{\Rightarrow}
\newcommand{\la}{\leftarrow}
\newcommand{\lla}{\longleftarrow}
\newcommand{\lra}{\longrightarrow}
\newcommand{\hra}{\hookrightarrow}

\newcommand{\rlas}{\rightleftarrows}

\newcommand{\CC}{\mathcal{C}}
\newcommand{\DD}{\mathcal{D}}

\newcommand{\SSc}{\mathcal{S}_{\bullet}}

\newcommand{\Rfd}{\mathcal{R}^{\it fd}}

\newcommand{\Ch}{{\mathcal C}h}
\newcommand{\Chfd}{{\mathcal C}h^{\it fd}}

\newcommand{\colim}[1][]{\ifmmode\ifinner {\operatorname{colim}_{#1}}\,\else \underset{#1}{\operatorname{colim}}\, \fi\fi}
\newcommand{\hocolim}[1][]{\ifmmode\ifinner {\operatorname{hocolim}_{#1}}\,\else \underset{#1}{\operatorname{hocolim}}\,\fi\fi}
\newcommand{\holim}[1][]{\ifmmode\ifinner {\operatorname{holim}_{#1}}\,\else \underset{#1}{\operatorname{holim}}\,\fi\fi}

\newcommand{\wh}[1][{h}]{{\rm Wh}^{\mathbb Q}_{#1}} 
\newcommand{\bC}[1][\ast]{C'_{#1}}

\newcommand{\taus}{\ensuremath\tau^{s}}
\newcommand{\ts}{\ensuremath t^{s}}

\newcommand{\hofib}{\mathrm{hofib} }

\newcommand{\id}{{\ensuremath{\rm id}}}

\newcommand{\End}{\operatorname{End}}

\newcommand{\coker}{\operatorname{coker}}
\newcommand{\Cyl}{\operatorname{Cyl}}
\newcommand{\Cone}{\operatorname{Cone}}
\newcommand{\SB}{\mathcal{S}(B)}

\newcommand{\R}{\mathbb{R}}
\newcommand{\Q}{\mathbb{Q}}
\newcommand{\Z}{\mathbb{Z}}

\newcommand{\benu}{\begin{enumerate}}
\newcommand{\eenu}{\end{enumerate}}

\title{Higher torsion and secondary transfer of unipotent bundles}

\author[B. Badzioch]{Bernard Badzioch}
\address[]{Department of Mathematics, University at Buffalo, SUNY, Buffalo, NY}
\email{badzioch@buffalo.edu}

\author[W. Dorabia{\l}a]{Wojciech Dorabia{\l}a}
\address[]{Department of Mathematics, Penn State Altoona, Altoona, PA }
\email{wud2@psu.edu}

\begin{document}

\begin{abstract}
\noindent Given a unipotent bundle of smooth manifolds we construct its secondary transfer map 
and show  that this map determines the higher smooth torsion of the bundle. 
This approach to higher torsion provides a new perspective on some of its properties. 
In particular it yields in a natural way a formula for torsion of a composition of two bundles. 
\end{abstract}

\date{\bf 2015.04.26}

\maketitle

%%% TOC
\tableofcontents
%%%

%%%%%%%%%%%%%%%%%%%%%%%%%%%%
% INTRO
%%%%%%%%%%%%%%%%%%%%%%%%%%%%

\section{Introduction}
\label{INTRO SEC}

By a smooth bundle of manifolds we will understand here a smooth submersion $p\colon E\to B$
where $E$ and $B$ are smooth compact manifolds. A smooth bundle with the fiber $F$ 
is unipotent if  $B$ is path connected and the graded vector space $H_{\ast}(F; \Q)$ admits 
a filtration such that  $\pi_{1}B$  acts trivially on the filtration quotients. Igusa and Klein \cite{IgusaBook, KleinMorse}  showed using fiberwise Morse theory
 that to any unipotent bundle one can associate the higher torsion invariant which depends not only on 
the topological structure of $p$, but also on its smooth structure.  Higher torsion proved to be a useful 
tool in the study of smooth bundles. In \cite{IguICM} Igusa showed, for example, that it can be used 
to detect exotic disc bundles constructed by Hatcher.

 In \cite{BDW} and \cite{BDKW} the present authors in collaboration 
 with Klein and Williams extended ideas of Dwyer, Weiss, and Williams \cite{DWW} 
 to obtain an alternative construction of torsion
 of unipotent bundles  based on the machinery of homotopy theory.
This construction can be briefly described as follows. Let $K(\Q)$ be the infinite 
loop space underlying the 
algebraic $K$-theory spectrum of the field of rational numbers. Given a smooth bundle 
$p\colon E\to B$
we can construct a map $c_{p}\colon B\to K(\Q)$ which, roughly speaking, assigns to each $b\in B$ 
the point of $K(\Q)$ represented by the singular chain complex $C_{\ast}(p^{-1}(b); \Q)$. 
The smooth Riemann-Roch theorem 
of \cite{DWW} implies that $c_{\rho}$ admits a factorization
\begin{equation*}
\label{DWW DIAG}
\begin{tikzpicture}
\matrix (m) 
[matrix of math nodes, row sep=3em, column sep=3em, text height=1.5ex, text depth=0.25ex]
{
& Q(E_{+}) \\
B & K(\Q) \\ 
};
\path[->, thick, font=\scriptsize]
(m-2-1) 
edge node[auto] {$p^{!}$} (m-1-2)
edge node[below] {$c_{p}$} (m-2-2)
(m-1-2)
edge node[auto] {$\lambda_{E}$} (m-2-2)
; 
\end{tikzpicture}
\end{equation*}
where $Q(E_{+})= \Omega^{\infty}\Sigma^{\infty}(E_{+})$, $p^{!}$ is the Becker-Gottlieb transfer
and $\lambda_{E}$ is the linearization map (\ref{LINEARIZATION NN}).

If $p$ is a unipotent bundle then the map $c_{p}$ is homotopic via a preferred homotopy to a constant map.
As a consequence we obtain a lift of $p^{!}$ to the space $\wh[s](E)$ which is the homotopy fiber of 
$\lambda_{E}$:
\begin{equation*}
\label{DWWTOR DIAG}
\begin{tikzpicture}
\matrix (m) 
[matrix of math nodes, row sep=2em, column sep=4em, text height=1.5ex, text depth=0.25ex]
{
& \wh[s](E) \\
& Q(E_{+}) \\
B & K(\Q) \\ 
};
\path[->, thick, font=\scriptsize]
(m-3-1) 
edge node[auto] {$p^{!}$} (m-2-2)
edge node[below] {$c_{p}$} (m-3-2)
(m-2-2)
edge node[auto] {$\lambda_{E}$} (m-3-2)
(m-1-2)
edge node[auto] {} (m-2-2)
; 
\path[->, thick, dashed, font=\scriptsize]
(m-3-1)
edge [bend left = 15] node[auto] {$\taus(p)$} (m-1-2.south west)
;
\end{tikzpicture}
\end{equation*}

 The lift $\taus(p)$ is the smooth torsion of the bundle $p$.

The homotopy class of $\taus(p)$ is an invariant of the smooth structure of $p$ in the following sense. 
If $p'\colon E'\to B$ is another smooth bundle and $f\colon E'\to E$ is a smooth bundle map  
then $f$ induces a map
$$f_{\ast}\colon \wh[s](E')\to \wh[s](E)$$
The map $f_{\ast}\tau^{s}(p')$ need not be homotopic to $\taus(p)$ in general, but this property 
does hold provided that $f$ is  a fiberwise diffeomorphism of bundles. 

The map $\tau^{s}(p)$ gives rise to a certain cohomology class 
$$t^{s}(p)\in \bigoplus_{k>0}H^{4k}(B; \R)$$  
\cite[4.10]{BDKW} which we will call the cohomological torsion of the bundle $p$.

In \cite[Section 9]{IgusaAx} Igusa showed that the cohomological torsion of the composition $pq$ 
can be, in some cases, computed from the torsion of the bundles $p$ and  $q$. Namely, if $q$ is an oriented  linear 
sphere bundle then
\begin{equation}
\label{EQ-TRAX}
\ts(pq)=\chi(F_{q})\ts(p)+{\rm tr}^{E}_{B}(\ts(q))
\end{equation}
where $\chi(F_{q})\in \Z$ is the Euler characteristic of the fiber of $q$ and 
$${\rm tr}^{E}_{B}\colon H^{\ast}(E; \R)\to  H^{\ast}(B; \R)$$
is the transfer map associated to $p$. In \cite{IgusaAx} Igusa calls  the formula  \eqref{EQ-TRAX} 
the transfer axiom and shows that taken together with a few other properties  it  uniquely determines 
the cohomological torsion. 

Igusa's arguments  can be used to show that the formula \eqref{EQ-TRAX} 
holds under more general conditions on $p$ and $q$, e.g. if dimensions of fibers of these bundles 
have the same parity. In \cite[Thm 7.1]{BDKW} we verified that the same is true in the case when $p$ 
is an arbitrary unipotent bundle and $q$ satisfies the assumptions of the Leray-Hirsch isomorphism 
theorem. One of our goals  in this paper is to show that this formula holds in general:

\begin{theorem}
\label{THM_TRAX}
The formula \eqref{EQ-TRAX} holds for any unipotent bundles $p\colon E\to B$ and 
$q\colon  D \to E$.  
\end{theorem}

In order to prove Theorem \ref{THM_TRAX} we develop a new construction of smooth torsion
based on the notion of the secondary transfer of unipotent bundles.  
The starting point  for this construction is the following fact:

\begin{theorem}
\label{MAIN1 THM}
Given a smooth bundle of compact manifolds $p\colon E\to B$ with fiber $F_{p}$ 
consider the diagram 
\begin{equation}
\label{MAIN DIAG}
\begin{tikzpicture}[baseline=(current bounding box.center)]
\matrix (m) 
[matrix of math nodes, row sep=3em, column sep=3em, text height=1.5ex, text depth=0.25ex]
{
Q(B_{+}) & Q(E_{+}) \\
K(\Q) & K(\Q) \\ 
};
\path[->, thick, font=\scriptsize]
(m-1-1) 
edge node[auto] {$Q(p^{!})$} (m-1-2)
edge node[anchor=east] {$\lambda_{B}$} (m-2-1)
(m-2-1) 
edge node[below] {$\chi(F_{p})$} (m-2-2)
(m-1-2)
edge node[auto] {$\lambda_{E}$} (m-2-2)
; 
\end{tikzpicture}
\end{equation}
where the lower horizontal map is given by the multiplication by the Euler characteristic $\chi(F_{p})\in \Z$ of 
$F_{p}$ and the upper horizontal map is the Becker-Gottlieb transfer of $p$.  
If $p$ is a unipotent bundle then this diagram commutes up to 
a preferred homotopy 
$$\eta_{p}\colon Q(B_{+})\times [0,1]\to K(\Q)$$
\end{theorem}

As an intermediate step in the proof of this result it will be convenient to work in a more general 
setting of unipotent fibrations, i.e. fibrations $p\colon E \to B$ satisfying some finiteness assumptions 
and such that the action of $\pi_{1}(B)$ on homology of the fiber satisfies the same unipotency
condition as in the case of unipotent bundles (see Definition \ref{UNIPFIBR-DEF}). We show 
(\ref{MAIN1-PROP}) that for any unipotent fibration an analog Theorem \ref{MAIN1 THM} holds,  
with  spaces $Q(B_{+})$ and $Q(E_{+})$ replaced with Waldhausen's algebraic $K$-theory spaces 
$A(B)$ and $A(E)$, and with $A$-theory transfer taken in place of the Becker-Gottlieb transfer.

For a unipotent bundle $p\colon E\to B$ the homotopy $\eta_{p}$ defines a map of homotopy fibers
$$\wh[s](p^{!})\colon \wh[s](B)\to \wh[s](E)$$
We call this map the smooth secondary transfer of the bundle $p$.  Likewise, for any unipotent fibration 
$p$ we construct its homotopy secondary transfer 
$$\wh[h](p^{!})\colon \wh[h](B)\to \wh[h](E)$$
where $\wh[h](B) = \hofib( A(B)\to K(\Q))$. 

The smooth secondary transfer shares some of the basic properties of the Becker-Gottlieb transfer. 
It is are additive (\ref{SADD THM})  and it preserves composition of bundles:

\begin{theorem}
\label{MAIN2 THM}
If $p\colon E\to B$ and $q\colon D\to E$ are unipotent bundles then 
$$\wh[s]((pq)^{!})\simeq \wh[s](q^{!})\circ \wh[s](p^{!})$$
\end{theorem}

Analogous additivity and composition properties hold for the homotopy secondary transfer 
(\ref{HADD THM}, \ref{HMAIN2-THM}).

The relationship between the smooth secondary transfer and the smooth torsion is as follows. 
If $B$ is a compact, smooth manifold then the identity map $\id_{B}\colon B \to B$
is a unipotent bundle. We have  

\begin{theorem}
\label{MAIN3 THM}
If $p\colon E\to B$ is a unipotent bundle then 
\begin{equation*}
\tau^{s}(p)\simeq \wh[s](p^{!})\circ\tau^{s}(\id_{B})
\end{equation*}
\end{theorem}
This shows that  the smooth secondary transfer of a unipotent bundle determines 
the smooth torsion of  the bundle. Since smooth torsion can distinguish bundles 
that are fiberwise homotopy equivalent, but not fiberwise diffeomorphic Theorem \ref{MAIN3 THM}
implies that the smooth secondary transfer carries information about the smooth structure of 
a bundle. In Proposition \ref{HOMOT INV PROP} we show that, in contrast, 
the homotopy secondary transfer is invariant with respect to fiberwise homotopy equivalences 
of unipotent fibrations.

Combining Theorems \ref{MAIN3 THM} and \ref{MAIN2 THM} we obtain

\begin{corollary}
\label{MAIN3 COR}
If $p\colon E \to B$ and $q\colon D \to E$ are unipotent bundles then 
$$\tau^{s}(pq)\simeq \wh[s](q^{!})\circ\tau^{s}(p)$$
\end{corollary}
 Theorem \ref{THM_TRAX} can be obtained as a direct consequence of this fact. 
 Notice that  in this way  we exhibit  the simple principle underlying 
the formula \eqref{EQ-TRAX}:  the torsion of a composition of bundles $p$ and $q$ 
is a composition of two maps, one depending  on $p$ and the other on $q$.

\begin{note}
As it was pointed out to us by the referee the main results of this paper bear some resemblance 
to the work of Lott and Bunke on the secondary $K$-theory pushforward map. 
In \cite{LottSecondary} Lott constructed for a smooth manifold 
$B$ the secondary $K$-theory group $\bar{K}^{0}(B)$ which is generated by flat complex vector 
bundles over $B$ with trivial Borel classes. He also constructed for a smooth bundle 
$p\colon E\to B$ a pushforward map $p_{!}\colon \bar{K}^{0}(E) \to \bar{K}^{0}(B)$. 
Just as the secondary smooth transfer considered  in this paper contains information 
about the smooth torsion of a bundle, the construction of  the pushforward map involves 
higher analytic torsion forms of Bismut and Lott \cite{Bismut-Lott}. Lott's pushforward map 
was studied by Bunke \cite{Bunke-SecLott} who showed that  it preserves composition of bundles: 
if $q\colon D \to E$ and $p\colon E\to B$ are smooth bundles then $(pq)_{!} = p_{!}q_{!}$. 
This  parallels our Theorem \ref{MAIN2 THM}. Beside the difference in setting 
between Lott's and Bunke's results and the ones described in this paper  the direction of their work 
is opposite to ours. While Lott's construction of the pushforward map uses explicitly the 
analytic torsion form, we construct the secondary transfer first and then show that the 
smooth torsion of a bundle can be recovered from it. Also, while we  use the composition property 
of the secondary transfer to obtain the composition formula for  cohomological torsion (\ref{EQ-TRAX}), 
Bunke derives his result using a theorem of Ma \cite{Ma-CompTorsion} which states that an analog of  
the formula (\ref{EQ-TRAX}) holds for higher analytic torsion. 

\end{note}

\begin{org}
Section \ref{TECHNICAL  SEC} contains a brief review of Waldhausen categories
which provide the technical setting for the majority of  construction of this paper. In Section  
\ref{LINEARIZATION-SEC} we take a closer look at the statements of Theorem \ref{MAIN1 THM} and 
its analog for unipotent fibrations, Theorem \ref{MAIN1-PROP}. The proof of both of these facts
uses a theorem of Brown,\!\!\!
\footnote{
See also \cite{IgusaTwisted} for a nice description of the relationship of Brown's work to 
the higher torsion of Igusa-Klein. 
}
which states that the singular chain 
complex of the total space of a  fibration is quasi-isomorphic to a twisted tensor product of the 
chain complexes of the base and the fiber.  In Section \ref{TWISTED SEC} we 
give an overview of this result and describe some properties of Brown's quasi-isomorphism.  
In Section \ref{SECTRANSFER-SEC}  we complete proofs of  Theorems 
\ref{MAIN1 THM} and \ref{MAIN1-PROP},   
which lets us complete the construction of the  secondary transfers $\wh[s](p^{!})$ and  $\wh(p^{!})$.
In \S \ref{HOMOTOPYINV SEC} we show that the homotopy secondary transfer if 
a fiberwise homotopy invariant.   
In Section \ref{ADDITIVITY SEC} we obtain additivity formulas for both smooth and homotopy 
secondary transfers, and in Sections \ref{COMPUNIP SEC} and \ref{COMPOSITION  SEC} we 
study  secondary transfers of compositions of unipotent bundles and fibrations, which leads us to 
the proof of Theorem \ref{MAIN2 THM}. 
Finally, in Section \ref{TORSION SEC} we prove the 
relationship between the smooth secondary transfer and the smooth torsion of a bundle $p$
described by Theorem \ref{MAIN3 THM}. We also show that it implies  the statement on  
Theorem \ref{THM_TRAX}. 
Several arguments of the paper involve constructions of maps 
between homotopy fibers and constructions of homotopies of such maps.
The appendix  (\S \ref{HOMOTFIB APP}) gives the basic outline of such constructions. 
\end{org}
\textbf{Acknowledgement.} We would like to thank the referee for several useful comments that 
helped us improve the presentation of this paper.

%%%%%%%%%%%%%%%%%%%%%%%%%%%%
% TECHNICAL SEC
%%%%%%%%%%%%%%%%%%%%%%%%%%%%

\section{Technical setup}
\label{TECHNICAL  SEC}

A great majority of constructions described in this paper is set within the realm of Waldhausen 
categories \cite{Wal}, i.e. categories with distinguished classes of weak equivalences and cofibrations 
that satisfy certain axioms. Our basic setup in this respect will be largely the same  as that of 
\cite[Section 3]{BDKW}, so we summarize it here only briefly.  Given a Waldhausen category $\CC$ we will 
denote by $K(\CC)$ the $K$-theory of $\CC$. The  standard construction of  $K(\CC)$ proceeds using  
Waldhausen's $\SSc$--construction. For our purposes it will be more convenient though to use 
its variant, the $\SSc '$-construction described by Blumberg and Mandell in \cite[\S 2]{Blumberg-Mandel}.

\begin{nn}
\label{AKQ CONSTR NN}
We will work mainly with two specific  instances of Waldhausen categories.   
For a topological space $X$ the category $\Rfd(X)$ has as its objects homotopy finitely dominated 
retractive spaces over $X$, while its morphisms are maps of retractive spaces.  
It is a Waldhausen category with cofibrations given by closed embeddings having the homotopy 
extension property and weak 
equivalences defined as homotopy equivalences. The $K$-theory of  $\Rfd(X)$ is the 
Waldhausen algebraic $K$-theory of $X$ and it is denoted  by $A(X)$
\footnote
{If  $X$ is a path connected space then by abuse of notation by $\Rfd(X)$ we will understand the 
Waldhausen category 
of \emph{path connected} retractive spaces over $X$. From the perspective of $K$-theory this 
change is of little consequence: the functor that embeds this category into the category of all retractive 
spaces over $X$ induces a homotopy equivalence of the associated $K$-theory spaces.}.  

Next, by  $\Chfd(\Q)$  we will denote the category of homotopy finitely dominated chain complexes of $\Q$-vector spaces. This is a Waldhausen category with degreewise mono\-mor\-phisms as cofibrations and quasi-isomorphisms as weak equivalences. We will denote by $K(\Q)$ the $K$-theory of $\Chfd(\Q)$. 
This space describes the algebraic $K$-theory of the field of rational numbers:  
$K(\Q)\simeq \Omega BGL(\Q)^{+}$. 

In order to obtain a convenient combinatorial model of the space $Q(X_{+})$ we will use one more 
instance of a construction of $K$-theory which, while it does not come from a Waldhausen category,  
is very closely related to the $S_{\bullet}$-construction. We outline it briefly here and refer 
to Waldhausen's papers \cite{WALM1, WALM2} where it was originally developed and
\cite{BDW} for details. Let $X$ be a smooth 
compact manifold. A partition of $X\times I$ is a codimension $0$ submanifold  $P\subseteq X\times I$
such that $X \times [0, \frac{1}{3}]\subseteq P$ and that satisfies some further conditions. By $\mathcal{P}_{0}(X\times I)$
we will denote the poset of partitions ordered by inclusion. Similarly, for $k\geq 0$ we define $\mathcal{P}_{k}(X\times I)$ 
to be the poset of locally trivial  bundles of partitions over the standard simplex $\Delta^{k}$. These categories taken
together form a simplicial category $\mathcal{P}_{\bullet}(X\times I)$. 
In each category $\mathcal{P}_{k}(X\times I)$ we can introduce an analog of 
the Waldhausen category structure where, roughly, every morphism is a cofibration and weak equivalences
are the identity morphisms. While in general pushouts do not exists in $\mathcal{P}_{k}(X\times I)$ and so the 
$\mathcal S_{\bullet}$-construction cannot be performed in this setting, it is possible to use its 
variant, the $\mathcal{T}_{\bullet}$-construction to obtain a space  
$|\mathcal{T}_{\bullet}\mathcal{P}_{k}(X\times I)|$.   
 Denote by $|\mathcal{T}_{\bullet}\mathcal{P}_{\bullet}(X\times I)|$
the space obtained by applying the $\mathcal{T}_{\bullet}$-construction to each  
category $\mathcal{P}_{k}(X\times I)$ and then taking geometrical realization of the resulting simplicial space. 
The final step in the 
construction of $Q(X_{+})$ is stabilization. It is obtained by means of maps    
$|\mathcal{T}_{\bullet}\mathcal{P}_{\bullet}(X\times I^{m})|\to 
|\mathcal{T}_{\bullet}\mathcal{P}_{\bullet}(X\times I^{m+1})|$
that are induced by functors $\mathcal{P}_{k}(X \times I^{m}) \to \mathcal{P}_{k}(X \times I^{m+1})$
that, roughly, send a partition $P\subseteq X \times I^{m}$ to $P\times I \subseteq X \times I^{m+1}$. 
Waldhausen showed  that there exists a weak equivalence
$$Q(X_{+}) \simeq \Omega \ \colim[m] |\mathcal{T}_{\bullet}\mathcal{P}_{\bullet}(X\times I^{m})|$$
\end{nn}

\begin{nn}
A functor $F\colon \CC\to \DD$ of Waldhausen categories is exact if it preserves 
weak equivalences, cofibrations and pushouts of diagrams 
\begin{equation}
\label{PUSHOUT EQ}
c' \la c \ra c''
\end{equation}
where one of the morphisms is a cofibration. Any such functor induces a map of infinite loop spaces 
$K(F)\colon K(\CC)\to K(\DD)$. The advantage of working with the $\SSc'$-construction is that 
we can obtain the map $K(F)$ under more relaxed assumptions on the functor $F$. Namely, following the terminology of \cite[3.4]{BDKW} we will  say that 
a functor $F$ is almost exact if it preserves weak equivalences and cofibrations, and if it preserves 
pushouts of diagrams (\ref{PUSHOUT EQ}) up to a weak equivalence. An almost exact functor  
induces a functor of simplicial categories $F\colon \SSc'\CC \to \SSc'\DD$, and so  it yields a map $K(F)\colon K(\CC)\to K(\DD)$.
\end{nn}

\begin{nn}
As we have already mentioned  exact and almost exact functors between Waldhausen 
categories define maps between their associated $K$-theories. We will frequently need 
to construct homotopies of such maps. There are two main sources of such constructions. 
First, if $F, G\colon \CC \to \DD$ are (almost) exact functors then a natural weak equivalence 
$\varphi\colon F \Ra G$  defines a homotopy $K(\varphi)$ between the induced maps 
$K(F)$ and $K(G)$. Second, if $F_{i}\colon \CC \to \DD$ are (almost) exact functors for $i=0, 1, 2$
and  $\varphi\colon F_{0}\Ra F_{1}$, $\psi\colon F_{1}\Ra F_{2}$ are natural transformations 
such that 
\begin{equation*}
\label{COFFUNCT EQ}
F_{0}(c)\overset{\varphi}{\to} F_{1}(c) \overset{\psi}{\to} F_{2}(c)
\end{equation*}
is a cofibration sequence for each $c\in \CC$, then the additivity theorem of Waldhausen 
\cite[Theorem 1.4.2]{Wal} provides a homotopy $\mho\colon K(\CC)\times I \to K(\DD)$ 
between the map $K(F_{1})$ and 
$K(F_{0}\vee F_{2})$. A convenient combinatorial construction of this homotopy was given by
Grayson in \cite{Grayson}. 
\end{nn}

\begin{nn}
Our work will require us to go a step further beyond homotopies and consider homotopies of 
homotopies. If $f_{0}, f_{1}\colon X\to Y$ are two maps between topological spaces and
$h_{0}, h_{1}\colon X\times I\to Y$ are homotopies between $f_{0}$ and $f_{1}$ then by a 
homotopy of homotopies we will understand a map $H\colon X\times I\times I \to Y$ such that 
for each $t\in I$ the map $H_{t} = H(-, -, t)$ is a homotopy between $f_{0}$ and $f_{1}$ with 
$H_{0} = h_{0}$, $H_{1} = h_{1}$. 
In our constructions of homotopies of homotopies we will almost  always have  $Y = K(\Q)$. 
The homotopies of homotopies we will  consider will be obtained as an application of 
one of the following lemmas:

\begin{lemma}
\label{HOMOTOHOMOT1 LEMMA}
Let $\CC$ be a Waldhausen category,  let 
$F_{1}, F_{2}\colon \CC \to \Chfd(\Q)$ be almost exact functors,  and let  
$\varphi_{1}, \varphi_{2}\colon F_{1}\Ra F_{2}$ be natural weak equivalences. 
If  $\Phi$ is a natural chain homotopy between $\varphi_{1}$ and $\varphi_{2}$ then  $\Phi$
defines a homotopy of homotopies $K(\Phi)$ between $K(\varphi_{1})$ and $K(\varphi_{2})$.  
\end{lemma}

\begin{proof}
We start with the following observation. 
Assume that we are given three almost exact functors $G_{i}\colon \CC \to \Chfd(\Q)$, $i=0,1,2$ 
and two natural weak equivalences $\psi_{0}\colon G_{0} \Ra G_{1}$ and 
$\psi_{1}\colon G_{1}\Ra G_{2}$. In such situation we obtain a map 
$$K(\CC)\times \Delta^{2} \to K(\Q)$$
If we consider it as a family of maps $K(\CC)\to K(\Q)$ parametrized by $\Delta^{2}$
then vertices of $\Delta^{2}$ correspond to the functions $K(G_{i})$ and edges of 
$\Delta^{2}$ correspond to the homotopies $K(\psi_{0})$, $K(\psi_{1})$, and $K(\psi_{1}\psi_{0})$. 

For a chain complex $C$ let $\Cyl(C)$ denote the mapping cylinder of the identity function 
$\id \colon C \to C$ \cite[1.5.5]{Weibel}. We have $\Cyl(C)_{n} = C_{n}\oplus C_{n-1}\oplus C_{n}$. 
Let $j_{0}, j_{1}\colon C \to \Cyl(C)$ be the chain maps given by $j_{0}(c) = (c, 0, 0)$ and 
$j_{1}(c) = (0, 0, c)$. Recall that two chain maps $f_{0}, f_{1}\colon C \to D$ are chain homotopic if and only if
there exists a chain map $h\colon \Cyl(C) \to D$ such that $hj_{i} = f_{i}$ for $i=0, 1$. 

Going back to the setting of the lemma the natural chain homotopy  $\Phi$ gives a natural weak 
equivalence of functors $\overline{\Phi}\colon \Cyl(F_{1}) \Ra F_{2}$. Let  
$\overline{\Psi} \colon \Cyl(F_{1}) \Ra F_{1}$ denote the natural weak equivalence 
corresponding to the natural chain homotopy from the identity natural transformation 
$\id \colon F_{1}\Ra F_{1}$ to itself. We obtain the following commutative diagram:
\begin{equation*}
\begin{tikzpicture}[baseline=(current bounding box.center)]
\matrix (m) 
[matrix of math nodes, row sep=3em, column sep=3em, text height=1.5ex, text depth=0.25ex]
{
 & F_{1} &  \\
F_{1} & \Cyl(F_{1}) & F_{2} \\
 & F_{1} &  \\ 
};
\path[-implies,thick, font=\scriptsize]
(m-1-2) 
edge [double] node[anchor = west] {$j_{0}$} (m-2-2)
edge [double] node[anchor = south west] {$\varphi_{1}$} (m-2-3)
edge [double] node[anchor = south east] {$\id$} (m-2-1)
(m-2-2)
edge [double] node[anchor = south] {$\overline{\Phi}\ \ $} (m-2-3)
edge [double] node[anchor = south] {$\ \ \overline{\Psi}$} (m-2-1)
(m-3-2) 
edge [double] node[anchor = west] {$j_{1}$} (m-2-2)
edge [double] node[anchor= north west] {$\varphi_{2}$} (m-2-3)
edge [double] node[anchor= north east] {$\id$} (m-2-1)
; 
\end{tikzpicture}
\end{equation*}
Each vertex of this diagram corresponds to a functor $\CC \to \Chfd(\Q)$ and edges 
are natural weak equivalences of such functiors. By the observation above each 
commutative triangle in this diagram induces a map $K(\CC) \times \Delta^{2} \to K(\Q)$. 
Taken together these maps define a map $H\colon K(\CC) \times I^{2} \to K(\Q)$. 
Considering $H$ as a family of functions $K(\CC) \to K(\Q)$ parametrized by $I^{2}$ we obtain that
two adjacent edges of the square $I^{2}$ parametrize the homotopies $K(\varphi_{1})$ 
and $K(\varphi_{2})$ and each of  remaining two edges parametrizes the homotopy defined by the identity 
natural transformation $\id\colon F_{1} \Ra F_{1}$. Using the identity $\varphi_{i}\id = \varphi_{i}$ 
we can further modify $H$ to a homotopy of homotopies between $K(\varphi_{1})$ and $K(\varphi_{2})$.

\end{proof}

\begin{lemma}
\label{HOMOTOHOMOT2 LEMMA}
Assume that  $F_{i}, G_{i}\colon \CC \to \Chfd(\Q)$ ($i=0, 1, 2$) are almost exact functors,
and that we have a commutative diagram of natural transformations 
\begin{equation*}
\begin{tikzpicture}[baseline=(current bounding box.center)]
\matrix (m) 
[matrix of math nodes, row sep=2em, column sep=2em, text height=1.5ex, text depth=0.25ex]
{
F_{0}  & F_{1} & F_{2} \\
G_{0} & G_{1} & G_{2} \\
};
\path[-implies, thick, font=\scriptsize]
(m-1-1) 
edge [double] (m-1-2)
edge [double] node[anchor=east] {$\eta_{0}$} (m-2-1)
(m-1-2)
edge [double]  (m-1-3)
edge [double] node[auto] {$\eta_{1}$} (m-2-2)
(m-1-3)
edge [double] node[auto] {$\eta_{2}$} (m-2-3)
(m-2-1)
edge  [double] (m-2-2)
(m-2-2)
edge  [double] (m-2-3)
; 
\end{tikzpicture}
\end{equation*}
where both rows are cofibration sequences and vertical arrows are natural weak equivalences. 
This yields a diagram 
\begin{equation}
\label{HOMOTOHOMOT2 EQ}
\begin{tikzpicture}[baseline=(current bounding box.center)]
\matrix (m) 
[matrix of math nodes, row sep=2.5em, column sep=2.5em, text height=1.5ex, text depth=0.25ex]
{
K(F_{1})& K(F_{0} \vee F_{2}) \\
K(G_{1})& K(G_{0} \vee G_{2}) \\
};
\path[->,thick, font=\scriptsize]
(m-1-1) 
edge node[auto] {$\mho$} (m-1-2)
edge node[anchor=east] {$K(\eta_{1})$} (m-2-1)
(m-1-2)
edge node[auto] {$K(\eta_{0} \vee \eta_{2})$} (m-2-2)
(m-2-1)
edge node[auto] {$\mho$} (m-2-2)
; 
\end{tikzpicture}
\end{equation}
In this diagram every vertex represents a map $K(\CC)\to K(\Q)$ and each edge is 
a homotopy of such maps. In this setting there exists a homotopy of homotopies 
that fills this diagram, i.e. a homotopy of homotopies between the concatenation 
of $\mho$ with $K(\eta_{0} \vee \eta_{2})$ and the concatenation of $K(\eta_{1})$ with 
$\mho$.
\end{lemma}

\begin{proof}
Denote by $\mathcal{E}\Chfd(\Q)$ the Waldhausen category of short exact sequences in 
$\Chfd(\Q)$. The functors $F_{i}$ and $G_{i}$ define functors 
$F, G \colon \CC \to \mathcal{E}\Chfd(\Q)$, and the natural weak equivalences $\eta_{i}$ 
define a natural weak equivalence $\eta\colon F \Ra G$. On the level of $K$-theory this 
yields a homotopy 
$$K(\eta)\colon K(\CC) \times I \to K(\mathcal{E}\Chfd(\Q))$$
For $i=0, 1, 2$ let $ev_{i}\colon \mathcal{E}\Chfd(\Q) \to \Chfd(\Q)$ denote the functor given by 
$$ev_{i}(C_{0} \to C_{1} \to C_{2}) = C_{i}$$ 
The Waldhausen additivity theorem can be equivalently stated by saying that there exists a homotopy 
$$\mathcal{E}\mho\colon \colon K(\mathcal{E}\Chfd(\Q)) \times I \to K(\Q)$$
between the map $K(ev_{1})$ and $K(ev_{0}\vee ev_{2})$. The homotopy of homotopies in
the diagram (\ref{HOMOTOHOMOT2 EQ}) is obtained by composing the map 
$$K(\eta)\times \id_{I} \colon K(\CC)\times I \times I \to K(\mathcal{E}\Chfd(\Q))\times I$$
with  $\mathcal{E}\mho$.  
\end{proof}
\end{nn}

\begin{lemma}
Consider the following a commutative diagram of almost exact functors 
$\CC \to \Chfd(\Q)$ and their natural transformations:
\begin{equation}
\label{3X3-DIAG}
\begin{tikzpicture}[baseline=(current bounding box.center)]
\matrix (m) 
[matrix of math nodes, row sep=2em, column sep=2em, text height=1.5ex, text depth=0.25ex]
{
A_{0} & B_{0} & C_{0} \\
A_{1}& B_{1} &  C_{1} \\
A_{2} & B_{2} & C_{2} \\
};
\path[-implies, thick, font=\scriptsize]
(m-1-1) 
edge [double] (m-1-2)
edge [double]  (m-2-1)
(m-1-2)
edge [double]  (m-1-3)
edge [double] (m-2-2)
(m-1-3)
edge [double]  (m-2-3)
(m-2-1)
edge  [double] (m-2-2)
edge [double]  (m-3-1)
(m-2-2)
edge  [double] (m-2-3)
edge [double]  (m-3-2)
(m-2-3)
edge  [double] (m-3-3)
(m-3-1)
edge  [double] (m-3-2)
(m-3-2)
edge  [double] (m-3-3)
; 
\end{tikzpicture}
\end{equation}
Assume that each row and each column of this diagram is a short exact sequence of functors 
and  that the canonical natural transformation $ A_{1}\cup_{A_{0}} B_{0} \Ra B_{1}$
is a cofibration. In this situation we obtain two homotopies between the map $K(B_{1})$ and 
the map $K(A_{0})+K(C_{0}) + K(A_{2})+K(C_{2})$:
\begin{itemize}[leftmargin=*]
\item[1)] the homotopy $\mho_{rc}$  is obtained by applying the additivity theorem to the 
middle row which produces a homotopy $K(B_{1})\simeq K(A_{1}) + K(C_{1})$,  and then 
applying additivity theorem to the left and right columns which gives a homotopy 
$K(A_{1})+K(C_{1}) \simeq K(A_{0})+K(A_{2}) + K(C_{0})+K(C_{2})$. 
\item[2)] the homotopy $\mho_{cr}$ which is obtained in the same manner, but applying additivity 
to the middle column first, and then to the top and bottom rows. 
\end{itemize}
There exists a homotopy of homotopies between $\mho_{rc}$ and $\mho_{cr}$
\end{lemma}

\begin{proof}

Let $\mathcal{E}\Chfd(\Q)$ denote  the Waldhausen category of short exact sequences in 
$\Chfd(\Q)$ and let $\mathcal{E}_{2}\Chfd(\Q)$ be the Waldhausen category of short 
exact sequences in $\mathcal{E}\Chfd(\Q)$.  The diagram (\ref{3X3-DIAG}) 
can be interpreted as a functor  $F\colon\CC \to  \mathcal{E}_{2}\Chfd(\Q)$ while columns  
of this diagram define functors $A, B, C \colon \CC \to \mathcal{E}\Chfd(\Q)$. The additivity theorem 
applied to $\mathcal{E}_{2}(\Chfd(\Q))$ gives a homotopy 
$\mathcal{E}_{2}\mho \colon K( \mathcal{E}_{2}\Chfd(\Q)) \to K(\mathcal{E}\Chfd(\Q))$
such that the composition 
\begin{equation*}
\begin{tikzpicture}[baseline=(current bounding box.center)]
\matrix (m) 
[matrix of math nodes, row sep=2em, column sep=3em, text height=1.5ex, text depth=0.25ex]
{
K(\CC) \times I & K( \mathcal{E}_{2}\Chfd(\Q))\times I & K( \mathcal{E}\Chfd(\Q)) \\
};
\path[->,thick, font=\scriptsize]
(m-1-1) 
edge node[auto]  {$K(F)\times \id_{I}$} (m-1-2)
(m-1-2)
edge node[auto]  {$\mathcal{E}_{2}\mho$} (m-1-3)
; 
\end{tikzpicture}
\end{equation*}
if a homotopy between the map $K(B)$ and $K(A)+K(C)$. 
Let $\mathcal{E}\mho\colon K( \mathcal{E}\Chfd(\Q)) \times I \to K( \mathcal{E}\Chfd(\Q))$
denote the additivity homotopy described in the proof of Lemma \ref{HOMOTOHOMOT2 LEMMA}. 
The homotopy of homotopies between  $\mho_{rc}$ and $\mho_{cr}$ is then given by the composition 
\begin{equation*}
\begin{tikzpicture}[baseline=(current bounding box.center)]
\matrix (m) 
[matrix of math nodes, row sep=2em, column sep=4em, text height=1.5ex, text depth=0.25ex]
{
K(\CC) \times I^{2}
& K( \mathcal{E}_{2}\Chfd(\Q))\times I^{2}
&[-1em] K( \mathcal{E}\Chfd(\Q)) \times I  
&[-2em] K(\Q)\\
};
\path[->,thick, font=\scriptsize]
(m-1-1) 
edge node[auto]  {$K(F)\times \id_{I^{2}}$} (m-1-2)
(m-1-2)
edge node[auto]  {$\mathcal{E}_{2}\mho \times \id_{I} $} (m-1-3)
(m-1-3)
edge node[auto]  {$\mathcal{E}\mho$} (m-1-4)

; 
\end{tikzpicture}
\end{equation*}

\end{proof}

%%%%%%%%%%%%%%%%%%%%%%%%%%%%
% THE LINEARIZATION MAP
%%%%%%%%%%%%%%%%%%%%%%%%%%%%

\section{The linearization  and transfer maps}
\label{LINEARIZATION-SEC}

\begin{nn}
\label{LINEARIZATION NN}
In preparation for the proof of Theorem \ref{MAIN1 THM} we start this section by reviewing briefly  
the construction of the linearization map $\lambda_{B}\colon Q(B_{+})\to K(\Q)$. 
For an arbitrary space $B$ we have the assembly map $a_{B}\colon Q(B_{+})\to A(B)$.
The combinatorial construction of this map can be outlined as follows. Recall (\ref{AKQ CONSTR NN})
that we are working with a model of $Q(B_{+})$ built using the category of partitions, and that   
the space $A(B)$ is constructed using the category $\Rfd(B)$.
By definition a partition  $P \subseteq B\times I$ contains $B\times \{0\}$ as a subspace
which gives $P$ the structure of a retractive space over $B$. In this way we can regard $P$ as 
an object of  $\Rfd(B)$. 
It is possible to extend this assignment to all parametrized and stabilized partitions  
(using parametrized and stabilized retractive spaces over $B$ are an intermediate step). 
The map $a_{B}$ is induced by this assignment. We refer to \cite[\S 3]{BDW} for details.

Next, recall that the $K(\Q)$ was constructed from the Waldhausen category 
of chain complexes $\Chfd(\Q)$. Consider the functor 
$$\Lambda_{B}\colon \Rfd(B)\to \Chfd(\Q)$$
that assigns to a retractive space $X$ the relative singular chain complex $C_{\ast}(X, B)$. 
This functor is almost exact,   
so it induces a map $\lambda^{h}_{B}\colon A(B)\to K(\Q)$. We will call 
$\lambda^{h}_{B}$ the $A$-theory linearization. 
The linearization map 
$\lambda_{B}\colon Q(B_{+})\to K(\Q)$ is  given by the composition 
$$\lambda_{B} = \lambda^{h}_{B}a_{B}$$ 
\end{nn}

\begin{nn}
\label{TRANSFER CONSTR NN}
Assume  that we have a fibration $p\colon E\to B$. 
For a retractive space $X\in \Rfd(B)$ let $p^{\ast}X$ denote the pullback
$$p^{\ast}X\coloneqq \lim(X \ra B \overset{\  p}{\la} E)$$
The assignment $X\mapsto p^{\ast}X$ defines an exact functor
$\Rfd(B)\to \Rfd(E)$. We will call the induced map $A(p^{!})\colon  A(B)\to A(E)$ the $A$-theory 
transfer of $p$. 

If $p$ is a smooth bundle of manifolds then in a similar way we can define a map 
$Q(p^{!})\colon Q(B_{+}) \to Q(E_{+})$. This map is induced by functors of categories of partitions
$\mathcal{P}_{k}(B\times I^{m}) \to \mathcal{P}_{k}(E\times I^{m})$ that, roughly,  associate to a partition 
$P\subseteq B\times I^{m}$ the partition $(p\times \id)^{-1}(P) \subseteq E\times I^{m}$.
The map $Q(p^{!})$ obtained in this way coincides with the Becker Gottlieb transfer \cite[\S 4]{BDW}.

\end{nn}

\begin{nn} 
\label{MAIN1HOMOTRED-NN}
Let $p\colon E\to B$ be a fibration with a homotopy finitely dominated fiber $F_{p}$. 
The maps described above can be assembled into a diagram
\begin{equation}
\label{MAIN-DIAG2}
\begin{tikzpicture}[baseline=(current bounding box.center)]
\matrix (m) 
[matrix of math nodes, row sep=3em, column sep=3em, text height=1.5ex, text depth=0.25ex]
{
Q(B_{+}) & Q(E_{+}) \\
A(B) & A(E) \\
K(\Q) & K(\Q) \\ 
};
\path[->,thick, font=\scriptsize]
(m-1-1) 
edge node[auto] {$Q(p^{!})$} (m-1-2)
edge node[anchor=east] {$a_{B}$} (m-2-1)
(m-1-2)
edge node[auto] {$a_{E}$} (m-2-2)
(m-2-1)
edge node[auto] {$A(p^{!})$} (m-2-2)
edge node[anchor=east] {$\lambda^{h}_{B}$} (m-3-1)
(m-2-2)
edge node[anchor=west] {$\lambda^{h}_{E}$} (m-3-2)
(m-3-1) 
edge node[below] {$\chi(F_{p})$} (m-3-2)
; 
\end{tikzpicture}
\end{equation}

The outer rectangle in this diagram coincides  with the diagram  (\ref{MAIN DIAG}). 
If $p\colon E\to B$ is a smooth bundle then directly from the constructions described above 
it follows that the upper square commutes up to a 
homotopy induced by the natural transformation that for a partition $P\subseteq B\times I$
is given by the isomorphism 
$$(p\times \id)^{-1}(P) \to p^{\ast}P$$ 
of retractive spaces over $E$. We will denote this homotopy by $\mu_{p}$
%%%
\!\!\footnote{While all maps in the upper square of (\ref{MAIN-DIAG2}), i.e. the Becker-Gottlieb 
transfer, the $A$-theory transfer and assembly maps are defined for any fibration, smoothness of $p$ 
is essential for commutativity. This diagram does not commute, in general, when $p$ is a fibration. 
See e.g. the proof of Theorem F in \cite{KW}}.
%%%
\!\!\!.  In order to obtain Theorem \ref{MAIN1 THM} it is  then enough to show 
that the lower square of  (\ref{MAIN-DIAG2}) is homotopy commutative. We will 
show that this fact holds for any unipotent fibration $p$:
\end{nn}

\begin{definition}
\label{UNIPFIBR-DEF}
A fibration $p\colon E\to B$ is \emph{unipotent} if $B$ is a path connected space,  both $B$
and the fiber $F_{p}$ of $p$ have the homotopy type of a finite CW-complex, 
and $H_{\ast}(F_{p})$
admits a filtration by $\pi_{1}B$-modules such that  the action of $\pi_{1}B$ on the 
filtration quotients is trivial.  
\end{definition} 

\begin{theorem}
\label{MAIN1-PROP}
Let  $p\colon E\to B$ be a  unipotent fibration with fiber $F_{p}$. 
The diagram
\begin{equation}
\label{MAIN-DIAG3}
\begin{tikzpicture}[baseline=(current bounding box.center)]
\matrix (m) 
[matrix of math nodes, row sep=3em, column sep=3em, text height=1.5ex, text depth=0.25ex]
{
A(B) & A(E) \\
K(\Q) & K(\Q) \\ 
};
\path[->, thick, font=\scriptsize]
(m-1-1)
edge node[auto] {$A(p^{!})$} (m-1-2)
edge node[anchor=east] {$\lambda^{h}_{B}$} (m-2-1)
(m-1-2)
edge node[anchor=west] {$\lambda^{h}_{E}$} (m-2-2)
(m-2-1) 
edge node[below] {$\chi(F_{p})$} (m-2-2)
; 
\end{tikzpicture}
\end{equation}
commutes up to a preferred homotopy $\eta^{h}_{p}\colon A(B)\times I\to K(\Q)$.
\end{theorem}

In the next section we  describe some technical tools that we will use  in the proof
of this fact. The proof itself is given in Section \ref{SECTRANSFER-SEC}.

%%%%%%%%%%%%%%%%%%%%%%%%%%%%
% TWISTED SEC
%%%%%%%%%%%%%%%%%%%%%%%%%%%%

\ 

\section{Twisted tensor products}
\label{TWISTED SEC}

Consider the diagram (\ref{MAIN-DIAG3}). The lower horizontal map in this diagram 
can be described combinatorially as follows. Let $H_{\ast}(F_{p})$ be the chain 
complex of rational homology groups of $F_{p}$ with trivial differentials. The functor
$$- \otimes H_{\ast}(F_{p}) \colon \Chfd(\Q) \to \Chfd(\Q)$$
is exact and $\chi(F_{p})\colon  K(\Q) \to K(\Q)$ is the map induced by this functor.  
As a consequence the map $\chi(F_{p})\lambda^{h}_{B}$ comes from the functor 
$\Rfd(B) \to \Chfd(\Q)$ that associates to a retractive space $X$ the chain complex
$C_{\ast}(X, B)\otimes H_{\ast}(F_{p})$. On the other hand the composition $\lambda^{h}_{E}A(p^{!})$
is induced by the functor $\Rfd(B) \to \Chfd(\Q)$ that assigns to a space $X$ the chain 
complex $C_{\ast}(p^{\ast}X, E)$. The main ingredient of the proof of Theorem \ref{MAIN1-PROP} 
is  the fact that  under the assumption that $p\colon E\to B$ is a unipotent fibration for any 
$X\in \Rfd(B)$ we can construct a path in $K(\Q)$ joining the points represented  
by these two chain complexes. This path is natural to the extent that it gives rise to a homotopy filling the 
diagram (\ref{MAIN-DIAG3}). Our main tool in the construction of this path will be the  theorem 
of Brown \cite{Brown} which shows that  the chain complex of the total space of a fibration 
is quasi-isomorphic to a twisted tensor product of the chain complex of the base and the homology 
of the fiber. We begin this section by reviewing
the relevant notions in homological algebra. Subsequently we describe Brown's result and 
develop  some of its properties that we will need later on.

\begin{twisting}
Let  $A$ be a differential graded  $\Q$-algebra with multiplication $\mu\colon A\otimes A\to A$, 
and let $K$ be a d.g. $\Q$-coalgebra with comultiplication  $\nabla\colon K\to K\otimes K$. 
Given homomorphisms of graded vector spaces $\varphi, \psi \colon K\to A$
the cup product $\varphi\cup \psi\colon K\to A$ is given by the formula
$$\varphi\cup \psi := \mu(\varphi\otimes \psi)\nabla$$
If $M$ is a left  $A$-module with multiplication $\nu\colon A\otimes M\to M$
then for $\varphi$ as above and $c\in K\otimes M$ 
 the cap product $\varphi\cap c \in K\otimes M$ is given by 
$$\varphi\cap c := (\id_{K}\otimes \nu)(\id_{K}\otimes \varphi\otimes \id_{M})(\nabla\otimes {\id_{M}})(c)$$
For a fixed $\varphi$ the map 
$$\varphi \cap -\colon K\otimes M\to K\otimes M$$
is a homomophism of graded vector spaces. 
\end{twisting}

\begin{definition}
\label{TWISTING COCHAIN}
Let $A$, $K$, $M$ be respectively a d.g. $\Q$-algebra,  coalgebra,  and a left 
$A$-module as above.  

\vskip 2mm 

\noindent (i) A twisting cochain is a homomorphism of graded vector spaces
$\varphi\colon K\to A$  that lowers grading by $1$ and  satisfies the identity 
$$\partial \varphi -\varphi\partial +\varphi\cup \varphi =0$$

\noindent (ii) If $\varphi\colon K\to A$ is a twisting cochain then  the twisted tensor product 
$K\otimes_{\varphi}M$ is a chain complex such that $K\otimes_{\varphi}M = K\otimes M$ 
as a graded vector space, and the differential in $K\otimes_{\varphi}M$ is given by
\begin{equation}
\label{TWISTED DIFF}
\partial_{\varphi}:= \partial\otimes\id + \id\otimes \partial + \varphi\cap -
\end{equation}
\end{definition}

\ 

\begin{twistedchain}
Let $X$ be a  topological space with a basepoint $x_{0}$. By $\bC(X)$ we will denote the 
subcomplex of the singular chain complex  $C_{\ast}(X)$ with coefficients in $\Q$ 
generated by all singular
simplices $\sigma\colon \Delta^{n}\to X$ that send all vertices of $\Delta^{n}$ into $x_{0}$. 
If $X$ is a path connected space then the inclusion $\bC(X)\hra C_{\ast}(X)$ is a chain 
homotopy equivalence. The chain complex $\bC(X)$  can be equipped with the usual 
d.g. coalgebra structure with comultiplication $\nabla\colon \bC(X)\to \bC(X)\otimes \bC(X)$ defined
by 
$\nabla(\sigma) = \sum_{i=0}^{n}f_{i}(\sigma)\otimes b_{n-i}(\sigma)$
where $\sigma\in \bC[n](X)$ is a singular $n$-simplex and $f_{i}(\sigma)$, $b_{n-i}(\sigma)$ denote, 
respectively, the front $i$-th face and the back $(n-i)$-th face of $\sigma$.  

For a space $X$ we can also consider  its  associated d.g.  homology algebra
$\End(H_{\ast}(X))$ defined as follows.  Let $\End_{n}(H_{\ast}(X))$ denote 
the vector space of all maps of graded vector spaces $H_{\ast}(X)\to H_{\ast}(X)$ that increase
the grading by $n$, and let 
$$\End(H_{\ast}(X)) = \bigoplus_{n\geq 0} \End_{n}(H_{\ast}(X))$$
We view $\End(H_{\ast}(X))$ as a chain complex with trivial differentials. 
The d.g. algebra structure on  $\End(H_{\ast}(X))$ comes from composition 
of maps. Naturally $H_{\ast}(X)$ is a module over this d.g. algebra. 

The main result of \cite{Brown} says that  given a fibration $p\colon E\to B$
with a path connected base space and a fiber $F_{p}$ we can find a twisting cochain 
$\varphi_{p}\colon \bC(B)\to \End(H_{\ast}(F_{p}))$ such that the twisted 
tensor product $\bC(B)\otimes_{\varphi_{p}}H_{\ast}(F_{p})$ is naturally quasi-isomorphic to 
$C_{\ast}(E)$. For our purposes it will be convenient to state this fact in the following form. 
Let $\mathcal{S}_{\ast}$ denote the category of path connected, pointed spaces, 
and $\mathcal{S}_{\ast}\!\!\downarrow \! B$ be the over category of $\mathcal{S}_{\ast}$ over 
a space $B$. 
Given a fibration $p\colon E\to B$ and an object $X\in \mathcal{S}_{\ast}\!\!\downarrow \!  B$
denote by $p_{X}\colon p^{\ast}X\to X$ the fibration induced from $p$. 
\end{twistedchain}
 
\begin{theorem}
\label{BROWN-THM}
Let $p\colon E\to B$  be a fibration with a path connected, pointed  base space $B$ and  a fiber $F_{p}$. 

1) For every $X\in \mathcal{S}_{\ast}\!\downarrow \! B$ there exists a twisting cochain 
$$\varphi_{p_{X}}\colon \bC(X)\to \End(H_{\ast}(F_{p}))$$
and a quasi-isomorphism 
$$\beta_{p_{X}}\colon  C_{\ast}(p^{\ast}X) \overset{\simeq}{\lra} 
\bC(X)\otimes_{\varphi_{p_{X}}}H_{\ast}(F_{p})
$$

2) On $\bC[1](X)$ the twisting cochain 
$\varphi_{p_{X}} \colon \bC[1](X) \to \End_{0}(H_{\ast}(F_{p}))$ is given as follows. 
If $\sigma$ is a singular simplex in $\bC[1](X)$ then $\sigma$ is a loop in $X$, and so it represents 
an element $[\sigma]\in \pi_{1}X$.  For $z\in H_{\ast}(F_{p})$ we have 
$$\varphi_{p_{X}}(\sigma)(z) = [\sigma] z - z$$
where the product $[\sigma]z$ is defined by the action of $\pi_{1}X$ on $H_{\ast}(F_{p})$. 

3) The assignment $X\mapsto \bC(X)\otimes_{\varphi_{p_{X}}}H_{\ast}(F_{p})$ defines a 
functor 
$$F\colon \mathcal{S}_{\ast}\!\downarrow \! B \to \Ch(\Q)$$
where $\Ch(\Q)$ is the category of chain complexes over $\Q$. 
If $f\colon X\to Y$ is a morphism in $\mathcal{S}_{\ast}\!\downarrow \! B$ then 
$F(f) = f_{\ast}\otimes \id_{H_{\ast}(F_{p})}$. 

4) The quasi-isomorphisms $\beta_{p_{X}}$ define a natural transformation of functors. 
\end{theorem}

\begin{notation}
\label{TWIST-NOTATION}
For simplicity from now on we will write $\bC(X)\otimes_{p_{X}}H_{\ast}(F_{p})$ to denote 
the complex $\bC(X)\otimes_{\varphi_{p_{X}}} H_{\ast}(F_{p})$.
\end{notation}

\begin{nn}
\label{BROWN-NN}
While we refer to Brown's paper \cite{Brown} for the proof of Theorem \ref{BROWN-THM}, 
a few comments will be useful later on. Brown constructs the quasi-isomorphisms 
$\beta_{p_{X}}$ in two stages. 
First, he shows that given a path connected space $X$ one 
can construct a twisting cochain 
$\psi_{X}\colon \bC(X)\to C_{\ast}(\Omega X)$ \cite[Theorem 4.1]{Brown}, 
where the d.g. algebra structure on $C_{\ast}(\Omega X)$ is obtained by composing 
the Eilenberg-Zilber map and the map $C_{\ast}(\Omega X\times \Omega X) \to C_{\ast}(\Omega X)$ 
induced by the loop multiplication.  The twisting cochain $\psi_{X}$ depends 
on the space $X$ only, not on a fibration over $X$. Moreover, the construction of $\psi_{X}$
is natural on the category of path connected spaces. 
The action of $\Omega X$ on the fiber $F_{p}$ of $p_{X}$ defines a $C_{\ast}(\Omega X)$-module
structure on $C_{\ast}(F_{p})$. Brown shows \cite[Theorem 4.2]{Brown} that the twisted tensor product 
$\bC(X)\otimes_{\psi_{X}}C_{\ast}(F_{p})$ is chain homotopy equivalent to $C_{\ast}(p^{\ast}X)$
via a chain homotopy equivalence natural in $X$.
 A minor technical point here is that in order to get 
a suitable action of $\Omega X$ on $F_{p}$ one needs to specify a weakly transitive lifting function
for the fibration $p$. 
This can be taken care of by first replacing the fibration $p\colon E\to B$
by the homotopy equivalent fibration  $\tilde p \colon   PB\times_{B} E \to B$ where $PB$ is 
the space of Moore paths in $B$. 
The fibration $\tilde p$ admits a canonical lifting function \cite[p.225]{Brown} which  can 
be used to get a chain homotopy equivalence  
$$
C_{\ast}(p^{\ast}X)  \overset{\simeq}{\lra} C_{\ast}(\tilde p^{\ast} X) \overset{\simeq}{\lra} 
\bC(X)\otimes_{\psi_{X}}C_{\ast}(F_{\tilde p}) 
 $$
where $F_{\tilde p}$ is the fiber of $\tilde p$. Since $F_{\tilde p}\simeq F_{p}$ we have 
$C_{\ast}(F_{\tilde p})\simeq C_{\ast}(F_{p})$. 
%%% FOOTNOTE
\!\!\footnote{In \cite{Brown} Brown gives a quasi-isomorphism going in the opposite direction, 
$\bC(X)\otimes_{\psi_{X}}C_{\ast}(F_{\tilde p}) \overset{\simeq}{\to} C_{\ast}(\tilde p^{\ast} X)$. 
However, since his argument relies on the method of acyclic models it also produces a natural 
homotopy inverse of that map, and we work here with this inverse for convenience.}
%%%

To complete the construction of $\beta_{p_{X}}$ it suffices to show that we have quasi-iso\-mor\-phisms
$$
\bC(X)\otimes_{\psi_{X}}C_{\ast}(F_{\tilde p})
\overset{\simeq}{\lra
}
\bC(X)\otimes_{p_{X}}H_{\ast}(F_{p}) 
$$
One can proceed as follows. 
Using the fact that we deal here with chain complexes over a field we can find chain maps 
$$j_{F_{p}}\colon H_{\ast}(F_{p}) \rlas C_{\ast}(F_{\tilde p}) \colon r_{F_{p}}$$
such that for $z\in H_{\ast}(F_{p})\cong H_{\ast}(F_{\tilde p})$ the element  
$j_{F_{p}}(z)\in C_{\ast}(F_{\tilde{p}})$
is a chain representing $z$, $r_{F_{p}}j_{F_{p}}= \id_{H_{\ast}(F_{p})}$, and 
$j_{F_{p}}r_{F_{p}}\simeq \id_{C_{\ast}(F_{\tilde{p}})}$. 
The maps $j_{F_{p}}$ and $r_{F_{p}}$ define a strong deformation retraction of (untwisted) tensor products
$$
\\id\otimes j_{F_{p}}\colon \bC(X)\otimes H_{\ast}(F_{p}) \rlas 
\bC(X)\otimes C_{\ast}(F_{\tilde p}) \colon \id \otimes r_{F_{p}}
$$
The Basic Perturbation Lemma (see e.g. \cite[2.6]{Lambe-Stasheff})  shows that in such situation 
there is  a twisting cochain $\varphi_{p_{X}}\colon \bC(X)\to \End(H_{\ast}(F_{p}))$
and a strong deformation retraction of twisted tensor products
$$
(\id\otimes j_{F_{p}})^{\infty}\colon \bC(X)\otimes_{p_{X}} H_{\ast}(F_{p}) \rlas 
\bC(X)\otimes_{\psi_{X}} C_{\ast}(F_{\tilde p}) \colon (\id \otimes r_{F_{p}})^{\infty}
$$
Following our convention (\ref{TWIST-NOTATION}) by $\bC(X)\otimes_{p_{X}} H_{\ast}(F_{p})$
we denote here  the complex $ \bC(X)\otimes_{\varphi_{p_{X}}} H_{\ast}(F_{p})$.
We define $\beta_{p_{X}}$ as the composition 
\begin{equation*}
\begin{tikzpicture}[baseline=(current bounding box.center)]
\matrix (m) 
[matrix of math nodes, row sep=3em, column sep=2em, text height=1.5ex, text depth=0.25ex]
{
\beta_{p_{X}}\colon  C_{\ast}(p^{\ast}X) 
 &  \bC(X)\otimes_{\psi_{X}}C_{\ast}(F_{\tilde p}) 
 & [2em] \bC(X)\otimes_{p_{X}}H_{\ast}(F_{p}) \\
};
\path[->, thick, font=\scriptsize]
(m-1-1)
edge  (m-1-2)
(m-1-2)
edge node[auto] {$(\id\otimes r_{F_{p}})^{\infty}$} (m-1-3)
; 
\end{tikzpicture}
\end{equation*}
\end{nn}

\begin{homolfiltr}
\label{HOMOLFILTR-NN}
Let $\varphi\colon K\to A$ be a twisting cochain and let $M$ be an $A$-module. Directly 
from the definition of a twisted tensor product  it follows that the 
chain complex $K\otimes_{\varphi}M$ admits an increasing  filtration
$$ U_{0} \subset U_{1} \subset   \dots \subset K\otimes_{\varphi}M $$ 
where $U_{n} := (\bigoplus_{q\leq n}K_{q})\otimes_{\varphi}M$. In the case where $M$ has 
trivial differentials we have also a decreasing filtration
$$K\otimes_{\varphi}M = L_{0} \supset L_{1} \supset \dots$$
given by $L_{n} := K\otimes_{\varphi}(\bigoplus_{q\geq n}M_{q})$. Since we will consider 
this filtration in the situation where $M$ is the homology of some chain complex we will call 
it the  homological filtration of $K\otimes_{\varphi}M$. 
\end{homolfiltr}

Let $p\colon E\to B$ be a fibration with fiber $F_{p}$, let $X\in \mathcal{S}_{\ast}\!\downarrow \! B$,
and let $\{L_{n}(p_{X})\}$ denote the homological filtration of the chain complex
$\bC(X)\otimes_{p_{X}}H_{\ast}(F_{p})$. We will need an explicit description 
of the quotients $L_{n}(p_{X})/L_{n+1}(p_{X})$. 
On the level of  graded vector spaces we have isomorphisms 
$$L_{n}(p_{X})/L_{n+1}(p_{X}) \cong \bC(X)\otimes H_{n}(F_{p})$$
In order to describe the differential in the filtration quotients notice that the differential 
in $\bC(X)\otimes_{p_{X}}H_{\ast}(F_{p})$ is given by 
$$\partial(\sigma\otimes z) = \partial\sigma \otimes z + 
\sum_{i=0}^{n}f_{i}(\sigma)\otimes \varphi_{p_{X}}(b_{n-i}(\sigma))(z)$$
where $\sigma$ is a singular simplex in $\bC(X)$, $z\in H_{\ast}(F_{p})$ and $f_{i}(\sigma)$, 
$b_{n-i}(\sigma)$ are, respectively,  the $i$-th front face and the $(n-i)$-th back face of $\sigma$.
Using part 2) of Theorem \ref{BROWN-THM}  we get from here 
$$\partial(\sigma\otimes z) = \partial\sigma \otimes z + 
f_{n-1}(\sigma)\otimes ([b_{1}(\sigma)]z - z) \ \ \ \  (\text{mod\ } L_{n+1}(p_{X}))$$
As a consequence we obtain 
\begin{proposition}
\label{HFILTR-PROP} 
Let $p\colon E\to B$ be a fibration with a fiber $F_{p}$. For $X\in  \mathcal{S}_{\ast}\!\downarrow \! B$
let $\bC(X)\otimes_{\pi_{1}X}H_{n}(F_{p})$ denote the chain complex such that 
$$(\bC(X)\otimes_{\pi_{1}X}H_{n}(F_{p}))_{k}= \bC[k-n](X)\otimes H_{n}(F_{p})$$
with  differential  given by 
$$\partial(\sigma\otimes z) = \partial\sigma \otimes z + 
f_{n-1}(\sigma)\otimes ([b_{1}(\sigma)]z - z)$$
for a singular simplex $\sigma\in \bC(X)$ and $z\in H_{n}(F_{p})$. We have a canonical isomorphism
$$\bC(X)\otimes_{\pi_{1}X}H_{n}(F_{p})\cong L_{n}(p_{X})/L_{n+1}(p_{X})$$
\end{proposition}

\ 

\begin{mapsfibr}
Assume that that we have a map of fibrations over $B$:
\begin{equation*}
\begin{tikzpicture}[baseline=(current bounding box.center)]
\matrix (m) 
[matrix of math nodes, row sep=2.5em, column sep=2.5em, text height=1.5ex, text depth=0.25ex]
{
E & & D \\
& B & \\ 
};
\path[->, thick, font=\scriptsize]
(m-1-1)
edge node[auto] {$g$} (m-1-3)
edge node[anchor=north east] {$p$} (m-2-2)
(m-1-3)
edge node[anchor=north west] {$q$} (m-2-2)
; 
\end{tikzpicture}
\end{equation*}

For $X\in \mathcal{S}_{\ast}\!\downarrow \! B$ let  $g_{X}\colon p^{\ast}X \to q^{\ast}X$
be the map of the induced fibrations over $X$. Consider the diagram 
\begin{equation*}
\begin{tikzpicture}[baseline=(current bounding box.center)]
\matrix (m) 
[matrix of math nodes, row sep=3em, column sep=3em, text height=1.5ex, text depth=0.25ex]
{
C_{\ast}(p^{\ast}X) &  C_{\ast}(q^{\ast}X) \\
 \bC(X)\otimes_{p_{X}}H_{\ast}(F_{p}) &  \bC(X)\otimes_{q_{X}}H_{\ast}(F_{q})  \\ 
};
\path[->, thick, font=\scriptsize]
(m-1-1)
edge node[auto] {$g_{X\ast}$} (m-1-2)
edge node[anchor=east] {$\beta_{p_{X}}$} (m-2-1)
(m-1-2)
edge node[anchor=west] {$\beta_{q_{X}}$} (m-2-2)
;
\path[->, thick, densely dashed]  
(m-2-1) 
edge (m-2-2)
; 
\end{tikzpicture}
\end{equation*}
where $F_{p}$, $F_{q}$ denote, respectively, the fibers of $p$ and $q$. 
We would like to construct a  natural lower horizontal map such that the 
resulting  diagram commutes up to a homotopy. 
The obvious candidate for such map  is $\id\otimes (g|_{F_{p}})_{\ast}$, where the homomorphism
$(g|_{F_{p}})_{\ast}\colon H_{\ast}(F_{p})\to H_{\ast}(F_{q})$ is  induced by restriction of  
$g$ to the fibers, but this map  is not a chain map in general. We can, however, proceed  as follows. 
By the construction of quasi-isomorphisms $\beta_{p_{X}}$ and $\beta_{q_{X}}$
(\ref{BROWN-NN}) we have a diagram 
\begin{equation}
\label{GXINFTY-EQ}
\begin{tikzpicture}[baseline=(current bounding box.center)]
\matrix (m) 
[matrix of math nodes, row sep=3em, column sep=4em, text height=1.5ex, text depth=0.25ex]
{
C_{\ast}(p^{\ast}X) &  C_{\ast}(q^{\ast}X) \\
\bC(X)\otimes_{\psi_{X}}C_{\ast}(F_{\tilde{p}}) &  \bC(X)\otimes_{\psi_{X}} C_{\ast}(F_{\tilde{q}}) \\
 \bC(X)\otimes_{p_{X}}H_{\ast}(F_{p}) &  \bC(X)\otimes_{q_{X}}H_{\ast}(F_{q})  \\ 
};
\path[->, thick, font=\scriptsize]
(m-1-1)
edge node[auto] {$g_{X\ast}$} (m-1-2)
edge node[anchor=east] {$\simeq$} (m-2-1)
(m-1-2)
edge node[anchor=west] {$\simeq$} (m-2-2)
(m-2-1)
edge node[auto] {$\id\otimes (g|_{\tilde{F}_{p}})_{\ast}$} (m-2-2)
edge node[anchor=east] {$(\id\otimes r_{F_{p}})^{\infty}$} (m-3-1)
(m-2-2)
edge node[anchor=west] {$(\id\otimes r_{F_{q}})^{\infty}$} (m-3-2)
;
\path[->,  thick, font=\scriptsize, densely dashed]  
(m-3-1) 
edge node[auto] {$g_{X}^{\infty}$} (m-3-2)
; 
\end{tikzpicture}
\end{equation}

The compositions of the vertical maps give $\beta_{p_{X}}$ and $\beta_{q_{X}}$. The upper 
square commutes by the naturality properties of Brown's theorem \cite[Theorem 4.2]{Brown}. 
Recall that the map $(\id\otimes r_{F_{p}})_{\infty}$ is a part of the strong deformation retraction 
data 
$$(\id\otimes j_{F_{p}})^{\infty} 
\colon \bC(X)\otimes_{p_{X}} H_{\ast}(F_{p}) \rlas 
\bC(X)\otimes_{\psi_{X}} C_{\ast}(F_{\tilde p}) \colon
(\id\otimes r_{F_{p}})^{\infty}$$
Define a map $g_{X}^{\infty}\colon \bC(X)\otimes_{p_{X}}H_{\ast}(F_{p})
\to  \bC(X)\otimes_{q_{X}}H_{\ast}(F_{q})$ by
$$g_{X}^{\infty} \coloneqq 
(\id\otimes r_{Fq})^{\infty}\circ (\id\otimes (g|_{\tilde{F}_{p}})_{\ast}) \circ (id\otimes j_{F_{p}})^{\infty}$$ 
The lower square in the diagram (\ref{GXINFTY-EQ})
commutes then up to a chain homotopy. The Basic Perturbation Lemma gives explicit formulas
for this chain homotopy and for the maps $(\id\otimes r_{Fq})^{\infty}$ and $(id\otimes j_{F_{p}})^{\infty}$.
Direct computations involving these formulas yield the following fact:
\end{mapsfibr}
\begin{proposition}
\label{MAPSFIBR-PROP}
Let $B$ be a pointed, path connected space. Let $p\colon E\to B$ and $q\colon D\to B$
be fibrations, and let $g\colon E\to D$ be a map of fibrations. 

1) The maps $g_{X}^{\infty}$ define a natural transformation of functors 
$ \mathcal{S}_{\ast}\!\downarrow \! B \to \Ch(\Q)$. 

2) The diagram 
\begin{equation*}
\begin{tikzpicture}[baseline=(current bounding box.center)]
\matrix (m) 
[matrix of math nodes, row sep=3em, column sep=3em, text height=1.5ex, text depth=0.25ex]
{
C_{\ast}(p^{\ast}X) &  C_{\ast}(q^{\ast}X) \\
 \bC(X)\otimes_{p_{X}}H_{\ast}(F_{p}) &  \bC(X)\otimes_{q_{X}}H_{\ast}(F_{q})  \\ 
};
\path[->, thick, font=\scriptsize]
(m-1-1)
edge node[auto] {$g_{X\ast}$} (m-1-2)
edge node[anchor=east] {$\beta_{p_{X}}$} node[anchor=west]{$\simeq$} (m-2-1)
(m-1-2)
edge node[anchor=west] {$\beta_{q_{X}}$} node[anchor=east] {$\simeq$} (m-2-2)
(m-2-1) 
edge  node[anchor=south] {$g_{X}^{\infty}$} (m-2-2)
; 
\end{tikzpicture}
\end{equation*}
commutes up to a chain homotopy that is natural in $X$.

3) For $X\in \mathcal{S}_{\ast}\!\downarrow \! B$ consider the homological filtrations 
$\{L_{n}(p_{X})\}$ and $\{L_{n}(q_{X})\}$ of the complexes  
$\bC(X)\otimes_{p_{X}} H_{\ast}(F_{p})$ and, respectively, 
$\bC(X)\otimes_{q_{X}} H_{\ast}(F_{q})$  (\ref{HOMOLFILTR-NN}). 
The map $g_{X}^{\infty}$ preserves these filtrations. 
Moreover, for every $n$ the following diagram commutes:
\begin{equation*}
\begin{tikzpicture}[baseline=(current bounding box.center)]
\matrix (m) 
[matrix of math nodes, row sep=3em, column sep=3em, text height=1.5ex, text depth=0.25ex]
{
L_{n}(p_{X})/L_{n-1}(p_{X}) & L_{n}(q_{X})/L_{n-1}(q_{X}) \\
 \bC(X)\otimes_{\pi_{1}X}H_{n}(F_{p}) &  \bC(X)\otimes_{\pi_{1}X}H_{n}(F_{q})  \\ 
};

\path[->, thick, font=\scriptsize]
(m-1-1)
edge node[auto] {$g^{\infty}_{X}$} (m-1-2)
edge node[anchor=east] {$\cong$}  (m-2-1)
(m-1-2)
edge node[anchor=west] {$\cong$} (m-2-2)
(m-2-1) 
edge  node[anchor=south] {$\id\otimes (g|_{F_{p}})_{\ast}$} (m-2-2)
; 
\end{tikzpicture}
\end{equation*}
The vertical isomorphisms  in this diagram come from  Proposition \ref{HFILTR-PROP}. 
\end{proposition}

%%%%%%%%%%%%%%%%%%%%%%%%%%%%
%  SECTRANSFER SEC
%%%%%%%%%%%%%%%%%%%%%%%%%%%%

\section{The secondary transfer}
\label{SECTRANSFER-SEC}

We are now ready to give
\begin{proof}[Proof of Theorem \ref{MAIN1-PROP}]  \  
Consider the diagram (\ref{MAIN-DIAG3}). We want to construct a homotopy 
$\eta^{h}_{p}$ between the maps $\lambda^{h}_{E}A(p^{!})$ and $\chi(F_{p})\lambda^{h}_{B}$. 
Recall that  that the map $\lambda^{h}_{E}A(p^{!})$ is induced by the functor 
$$\Phi\colon \Rfd(B)\to \Chfd(\Q)$$
that assigns to a retractive space $X$ the relative chain complex $C_{\ast}(p^{\ast}X, E)$. 
while the map $\chi(F_{p})\lambda^{h}_{B}$ is induced by the functor
$$\Psi\colon \Rfd(B)\to \Chfd(\Q)$$
given by $\Psi(X) = C_{\ast}(X, B)\otimes H_{\ast}(F_{p})$.  
We will build a sequence of intermediate functors between $\Phi$ and $\Psi$ and 
show that maps induced by these functors can be connected by homotopies. 

First, for  $X\in \Rfd(B)$ the inclusion map $i_{X} \colon B\hra X$ induces an inclusion 
$\tilde \imath_{X}\colon E\hra p^{\ast}X$. 
The naturality of the quasi-isomorphisms $\beta_{p_{X}}$ described in Theorem \ref{BROWN-THM}
implies that we have a commutative diagram
\begin{equation*}
\begin{tikzpicture}[baseline=(current bounding box.center)]
\matrix (m) 
[matrix of math nodes, row sep=3em, column sep=4em, text height=1.5ex, text depth=0.25ex]
{
C_{\ast}(E) &  C_{\ast}(p^{\ast}X) \\
 \bC(B)\otimes_{p}H_{\ast}(F_{p}) &  \bC(X)\otimes_{p_{X}}H_{\ast}(F_{p})  \\
};
\path[->, thick, font=\scriptsize]
(m-1-1)
edge node[auto] {$\tilde \imath_{X\ast}$} (m-1-2)
edge node[anchor=east] {$\beta_{p}$} (m-2-1)
(m-1-2)
edge node[anchor=west] {$\beta_{p_{X}}$}  (m-2-2)
(m-2-1) 
edge  node[anchor=south] {$i_{X\ast}\otimes \id$} (m-2-2)
; 
\end{tikzpicture}
\end{equation*}
Define 
$$\bC(X, B)\otimes_{p_{X}}H_{\ast}(F_{p})\coloneqq \coker(i_{X\ast}\otimes \id)$$
The assignment $X\mapsto \bC(X, B)\otimes_{p_{X}}H_{\ast}(F_{p})$ defines an 
almost exact functor 
$$\Phi_{1}\colon \Rfd(B)\to \Chfd(\Q)$$ 
and the  natural quasi-isomorphisms 
$$\beta_{p_{X}}\colon C_{\ast}(p^{\ast}X, E) 
\overset{\simeq}{\lra} \bC(X, B)\otimes_{p_{X}}H_{\ast}(F_{p})
$$
define a natural weak equivalence
$\beta\colon \Phi \Ra \Phi_{1}$.
Denote by $K(\Phi_{1})\colon A(B)\to K(\Q)$ the  map induced by $\Phi_{1}$. The natural weak equivalence
$\beta$ defines a homotopy 
\begin{equation}
\label{1STHOMOT-EQ}
\lambda^{h}_{E}A(p^{!}) \simeq K(\Phi_{1})
\end{equation}
Next, since the map $i_{X\ast}\otimes\id$ preserves the homological filtration of  twisted 
tensor products we can define
\begin{equation*}
\begin{tikzpicture}
\matrix (m) 
[matrix of math nodes, row sep=3em, column sep=3em, text height=1.5ex, text depth=0.25ex]
{
L_{n}(p_{X}, p) \coloneqq \coker(L_{n}(p)  &  L_{n}(p_{X})) \\
};
\path[->, thick, font=\scriptsize]
(m-1-1)
edge node[auto] {$i_{X\ast}\otimes\id$} (m-1-2)
; 
\end{tikzpicture}
\end{equation*}
The chain complexes $L_{n}(p_{X}, p)$ form a decreasing filtration of $\bC(X, B)\otimes_{p_{X}}H_{\ast}(F_{p})$.
Proposition  \ref{HFILTR-PROP} shows that the filtration quotient $L_{n}(p_{X}, p)/L_{n+1}(p_{X}, p)$ 
can be identified with the chain complex 
$$\bC(X, B)\otimes_{\pi_{1}X}H_{n}(F_{p}) \coloneqq 
\coker (\bC(B)\otimes_{\pi_{1}B}H_{n}(F_{p}) \to \bC(X)\otimes_{\pi_{1}X}H_{n}(F_{p}))$$ 
Since $F_{p}$ is a homotopy finite space
we have $H_{q}(F_{p}) = 0$ for $q$ large enough, and so $\{L_{n}(p_{X}, p)\}$
is in fact a finite filtration.   The assignments $X\to L_{n}(p_{X}, p)$, and 
$X\mapsto \bC(X, B)\otimes_{\pi_{1}X}H_{n}(F_{p})$ define almost exact functors $\Rfd(B)\to \Chfd(\Q)$. 
These functors are connected by natural short exact sequences
\begin{equation}
\label{LNEXSEQ-EQ}
0\to L_{n+1}(p_{X}, p)\to L_{n}(p_{X}, p)\to  \bC(X, B)\otimes_{\pi_{1}X}H_{n}(F_{p}) \to 0
\end{equation}
Let $\Phi_{2}\colon \Rfd(B)\to \Chfd(\Q)$ denote the almost exact functor given by 
$$\Phi_{2}(X)  \coloneqq \bigoplus_{n}  \bC(X, B)\otimes_{\pi_{1}X}H_{n}(F_{p})$$
and let $K(\Phi_{2})\colon A(B)\to K(\Q)$ be the map induced by $\Phi_{2}$. Applying repeatedly 
Waldhausen's additivity theorem to the sequences (\ref{LNEXSEQ-EQ}) we obtain a
homotopy 
\begin{equation}
\label{2NDHOMOT-EQ}
K(\Phi_{1}) \simeq K(\Phi_{2})
\end{equation}

Assume now for a moment that $p\colon E\to B$ is a fibration with the trivial action of $\pi_{1}B$ 
on $H_{\ast}(F)$. In this case the action of $\pi_{1}X$ on $H_{\ast}(F_{p})$ is trivial as well, 
so we have isomophisms 
$$\Phi_{2}(X)\cong \bC(X, B)\otimes H_{\ast}(F_{p})$$
Since by assumption $X$ and $B$ are path connected spaces we also have natural 
quasi-isomorphisms
$$\bC(X, B)\otimes H_{\ast}(F_{p}) \simeq C_{\ast}(X, B)\otimes H_{\ast}(F_{p})= \Psi(X)$$
As a consequence for every $X\in \Rfd(B)$ we obtain $\Phi_{2}(X)\simeq \Psi(X)$ which induces 
a homotopy 
\begin{equation}
\label{3RDHOMOT-EQ}
K(\Phi_{2})\simeq \chi(F_{p})\lambda^{h}_{B}
\end{equation}
Concatenating the homotopies (\ref{1STHOMOT-EQ}), (\ref{2NDHOMOT-EQ}), and 
(\ref{3RDHOMOT-EQ}) we get the desired homotopy $\eta^{h}_{p}$.
 
If $p$ is an arbitrary unipotent fibration we need an additional step to pass between the maps 
$K(\Phi_{2})$ and $\chi(F_{p})\lambda^{h}_{B}$.  
In this  case the action of $\pi_{1}B$ need not be trivial, but $H_{\ast}(F_{p})$ admits a decreasing
filtration  $\{V^{i} \}$  such that  $V^{i}$  is a $\pi_{1}B$-module and
 the action of   $\pi_{1}B$ on the quotients $V^{i+1}/V^{i}$ is trivial.  
This defines a filtration $\{\bC(X, B)\otimes_{\pi_{1}X}V^{i}\}$ of the complex 
$\Phi_{2}(X)$. The quotients of this filtration are the (untwisted)
tensor products $\bC(X, B)\otimes (V^{i}/V^{i-1})$.
Define a functor $\Phi_{3}\colon \Rfd(B)\to \Chfd(\Q)$ by 
$$\Phi_{3}(X) \coloneqq \bigoplus_{i} \bC(X, B)\otimes (V^{i}/V^{i-1})$$
Naturally we also have a filtration  $\{\bC(X, B)\otimes V^{i}\}$  of the untwisted
tensor product $\bC(X, B)\otimes H_{\ast}(F_{p})$ and  $\Phi_{3}(X)$ is the direct sum 
of the quotients of this filtration. This means that using Waldhausen's additivity theorem 
(and the quasi-isomorphisms $\bC(X, B)\otimes H_{\ast}(F_{p})\simeq C_{\ast}(X, B)\otimes H_{\ast}(F_{p})$)
we get homotopies
\begin{equation}
\label{4THHOMOT-EQ}
K(\Phi_{2})\simeq K(\Phi_{3})\simeq  \chi(F_{p})\lambda^{h}_{B}
\end{equation}
The homotopy $\eta^{h}_{p}$ is then obtained as a concatenation of the homotopies 
(\ref{1STHOMOT-EQ}), (\ref{2NDHOMOT-EQ}), and 
(\ref{4THHOMOT-EQ}). 
\end{proof}

\begin{proof}[Proof of Theorem \ref{MAIN1 THM}] 
The homotopy $\eta_{p}$ is obtained by concatenating the homotopy $\eta^{h}_{p}$ and the homotopy 
$\mu_{p}$ described in (\ref{MAIN1HOMOTRED-NN}). 
\end{proof}

\

Let $C\in \Chfd(\Q)$ be a chain complex. Notice that  our construction of $K(\Q)$ lets us
identify $C$ with a point of $K(\Q)$.

\begin{definition}
\label{QWH-DEF}
Let $B$ be a path connected space and let $C\in \Chfd(\Q)$. Denote by  
$\wh[s](B)_{C}$ the homotopy fiber of the linearization map taken over the point 
$C\in K(\Q)$. 
$$\wh[s](B)_{C} \coloneqq \hofib(\lambda_{B}\colon Q(B_{+})\to K(\Q))_{C}$$
For the zero chain complex $0\in \Chfd(\Q)$ we will write $\wh[s](B)$ 
to denote $\wh[s](B)_{0}$.
\end{definition}

Let $p\colon E\to B$ be a unipotent bundle with a fiber $F_{p}$. 
By Theorem \ref{MAIN1 THM} for any $C\in \Chfd(\Q)$ we have a map 
$$\wh[s](B)_{C}\lra \wh[s](E)_{C\otimes H_{\ast}(F_{p})}$$
This gives rise to the following

\begin{definition}
\label{QSEC TRANSFER DEF}
The smooth secondary transfer of  a unipotent bundle  $p\colon E\to B$ 
is the map 
$$\wh[s](p^{!})\colon \wh[s](B)\lra \wh[s](E)$$
determined by the Becker-Gottlieb transfer $Q(p^{!})$ and the homotopy $\eta_{p}$
given by Theorem \ref{MAIN1 THM}.
\end{definition}

It will be convenient to consider a variant of this definition in the setting of unipotent 
fibrations:

\begin{definition}
\label{ASEC TRANSFER DEF}
For a path connected space $B$ let 
$$\wh(B) \coloneqq \hofib(\lambda^{h}_{B}\colon A(B)\to K(\Q))_{0}$$
The homotopy secondary transfer of  a unipotent fibration  $p\colon E\to B$ 
is the map 
$$\wh(p^{!})\colon \wh(B)\lra \wh(E)$$
determined by the transfer $A(p^{!})$ and the homotopy $\eta^{h}_{p}$ given by 
Theorem \ref{MAIN1-PROP}.
\end{definition}

\begin{note}
Let $p\colon E\to B$ be a unipotent fibration with a fiber $F_{p}$. 
The construction of the homotopy $\eta^{h}_{p}$ described in the proof of Theorem 
\ref{MAIN1-PROP} makes use of a choice of a strong deformation retraction 
$H_{\ast}(F_{p})\rlas C_{\ast}(F_{p})$
and a choice of a unipotent filtration $\{V^{i}\}$ of $H_{\ast}(F_{p})$. One can check though 
that the homotopy class of the map $\wh(p^{!})$ is independent of these choices, and 
so it depends on the fibration $p$ only. 
Likewise, if $p$ is a unipotent bundle then the homotopy class of the smooth secondary 
transfer $\wh[s](p^{!})$ depends only on the bundle $p$. 

\end{note}

%%%%%%%%%%%%%%%%%%%%%%%%%%%%
% HOMOTOPY INVARIANCE
%%%%%%%%%%%%%%%%%%%%%%%%%%%%

\section{Homotopy invariance of $\wh(p^{!})$}
\label{HOMOTOPYINV SEC}

Let   $f\colon E_{1}\to E_{2}$ be a map of topological spaces. Such map defines 
an exact functor of Waldhausen categories 
$$f_{\ast}\colon \Rfd(E_{1})\to \Rfd(E_{2})$$
given by $f_{\ast}(X) = X\cup_{E_{1}}E_{2}$ for $X\in \Rfd(E_{1})$. This functor in turn 
induces a map $f_{\ast}\colon A(E_{1}) \to A(E_{2})$.  Consider the diagram 
\begin{equation*}
\begin{tikzpicture}[baseline=(current bounding box.center)]
\matrix (m) 
[matrix of math nodes, row sep=3em, column sep=1.5em, text height=1.5ex, text depth=0.25ex]
{
A(E_{1})&  & A(E_{2}) \\
& K(\Q) &   \\
};
\path[->, thick, font=\scriptsize]
(m-1-1)
edge node[auto] {$f_{\ast}$} (m-1-3)
edge node[anchor=north east] {$\lambda^{h}_{E_{1}}$} (m-2-2)
(m-1-3)
edge node[anchor=north west] {$\lambda^{h}_{E_{2}}$}  (m-2-2)
; 
\end{tikzpicture}
\end{equation*}
Recall that the map $\lambda^{h}_{E_{1}}$ is induced by the functor $\Rfd(E_{1})\to \Chfd(\Q)$
given by $X\mapsto C_{\ast}(X, E_{1})$. Similarly,  the map $\lambda^{h}_{E_{2}}f_{\ast}$
comes from the functor defined by $X\mapsto C_{\ast}( X\cup_{E_{1}}E_{2}, E_{2})$.  The 
natural quasi-isomorphisms $C_{\ast}(X, E_{1})\to C_{\ast}( X\cup_{E_{1}}E_{2}, E_{2})$ 
define a homotopy $h_{f}$ between $\lambda^{h}_{E_{1}}$ and $\lambda^{h}_{E_{2}}f_{\ast}$. As 
a result we obtain a map 
$$f_{\ast}\colon \wh(E_{1})\lra \wh(E_{2})$$
Our goal in this section is to prove the following

\begin{proposition}
\label{HOMOT INV PROP}
For $i=1, 2$ let $p_{i}\colon E_{i}\to B$ be a unipotent fibration, and let $f\colon E_{1}\to E_{2}$
be a fiberwise homotopy equivalence. There is a homotopy 
$$f_{\ast}\wh(p_{1}^{!})\simeq \wh(p_{2}^{!})$$ 
\end{proposition}

\begin{proof}
Let $F_{p_{i}}$ denote the fiber of $p_{i}$. 
Using (\ref{HOMOTFIBHOMOT NN}) we see that in order to obtain the desired homotopy 
it is enough to construct the following data:

1) a homotopy $H^{A}_{f}\colon A(B)\times I\to A(E_{2})$ 
between $f_{\ast}A(p_{1}^{!})$ and $A(p_{2}^{!})$; 

2)  a homotopy $H^{K}_{f}\colon K(\Q)\times I \to K(\Q)$ between 
the maps $\chi(F_{p_{1}})$ and $\chi(F_{p_{2}})$; 

3) a homotopy of homotopies that fills the following diagram:
\begin{equation}
\label{HOMOTINV_MAPS}
\begin{tikzpicture}[baseline=(current bounding box.center)]
\matrix (m) 
[matrix of math nodes, row sep=4em, column sep=5em, text height=1.5ex, text depth=0.25ex]
{
\lambda^{h}_{E_{2}}f_{\ast} A(p_{1}^{!}) & \lambda^{h}_{E_{2}} A(p_{2}^{!}) \\
\lambda^{h}_{E_{1}} A(p_{1}^{!}) & \\
\chi(F_{p_{1}})\lambda^{h}_{B}  &  \chi(F_{p_{2}})\lambda^{h}_{B} \\
};
\path[solid, thick, font=\scriptsize]
(m-1-1)
edge node[auto] {$\lambda^{h}_{E_{2}}H^{A}_{f}$} (m-1-2)
edge node[anchor=east] {$h_{f}\circ(A(p_{1}^{!})\times \id_{I})$} (m-2-1)
(m-1-2)
edge node[anchor=west] {$\eta^{h}_{p_{2}}$} (m-3-2)
(m-2-1)
edge node[anchor=east] {$\eta^{h}_{p_{1}}$}  (m-3-1)

(m-3-1) 
edge  node[anchor=north] {$H^{K}_{f}(\lambda^{h}_{B}\times \id_{I})$} (m-3-2)
; 
\end{tikzpicture}
\end{equation}

Each vertex of this diagram  represents a map $A(B)\to K(\Q)$ and edges represent homotopies 
of such maps. 

1) \emph{Construction of $H^{A}_{f}$.} The map $f_{\ast}A(p_{1}^{!})$ comes from the 
functor $\Rfd(B)\to \Rfd(E_{2})$ given by $X\mapsto f_{\ast}p_{1}^{\ast}X$ while $A(p_{2}^{!})$
is induced by the functor $X\mapsto p_{2}^{\ast}X$.  Since $f$ is a fiberwise homotopy equivalence 
the  natural maps  $f_{\ast}p_{1}^{\ast}X \to p_{2}^{\ast}X$ induced by $f$ are weak equivalences, 
and so they define the homotopy $H^{A}_{f}$. 

2) \emph{Construction of $H^{K}_{f}$.} Recall that for $i=1, 2$ the map $\chi(F_{i})$ 
is induced by functor $\Chfd(\Q)\to \Chfd(\Q)$ given  by $C\mapsto C\otimes H_{\ast}(F_{i})$.
Since the map  $f|_{F_{1}}\colon F_{1}\to F_{2}$ is a homotopy equivalence it  induces 
an isomorphism of homology groups of the fibers 
 $$(f|_{F_{1}})_{\ast}: H_{\ast}(F_{1})\overset{\cong}{\lra} H_{\ast}(F_{2})$$
This gives a natural isomorphism of functors 
$$-\otimes H_{\ast}(F_{1}) \Ra - \otimes H_{\ast}(F_{2})$$
The homotopy $H_{f}^{K}$ is defined by this natural isomorphism. 

3) \emph{Construction of the homotopy of homotopies.}  Consider the following diagram:
\begin{equation*}
\begin{tikzpicture}[baseline=(current bounding box.center)]
\matrix (m) 
[matrix of math nodes, row sep=6em, column sep=6em, text height=1.5ex, text depth=0.25ex]
{
C_{\ast}(f_{\ast}p_{1}^{\ast}X, E_{2}) &  C_{\ast}(p_{2}^{\ast}X, E_{2}) \\
C_{\ast}(p_{1}^{\ast}X, E_{1}) & \\[-2em]
\bC(X, B)\otimes_{p_{1X}}H_{\ast}(F_{p_{1}}) &  \bC(X, B)\otimes_{p_{2X}}H_{\ast}(F_{p_{2}})  \\[-2em]  
 \bC(X, B)\otimes H_{\ast}(F_{p_{1}}) &  \bC(X, B)\otimes H_{\ast}(F_{p_{2}})  \\ 
};

\path[->, thick, font=\scriptsize]
(m-1-1)
edge node[auto] {} (m-1-2)
(m-1-2)
edge node[anchor=west] {$\beta_{p_{2X}}$} (m-3-2)
(m-2-1)
edge node[anchor=north ] {$f_{X\ast}$} (m-1-2)
edge node[anchor= east] {$\beta_{p_{1X}}$} (m-3-1)
edge  (m-1-1)
(m-3-1)
edge node[anchor= south] {$f_{X}^{\infty}$} (m-3-2)
(m-4-1)
edge node[anchor= north] {$\id\otimes (f|_{F_{1}})_{\ast}$} (m-4-2)
;
\path[solid,  thick, font=\scriptsize]  
(m-3-1) 
edge node[anchor=east] {additivity} (m-4-1)
(m-3-2) 
edge node[anchor=west] {additivity} (m-4-2)
; 

\node[font=\Large] at (0,2.7) {\ding{192}};
\node[font=\Large] at (0,0) {\ding{193}};
\node[font=\Large] at (0,-2.4) {\ding{194}};
\end{tikzpicture}
\end{equation*}

Each vertex of this diagram represents a functor $\Rfd(B)\to \Chfd(\Q)$. The edges represent 
natural weak equivalences, with the exception of the lowest vertical edges where the passage 
between functors is obtained using additivity.  The maps $\beta_{p_{1X}}$ and 
$\beta_{p_{2X}}$ are the Brown quasi-isomorpshisms (\ref{BROWN-THM}) and the maps 
$f_{X}^{\infty}$ come from Proposition \ref{MAPSFIBR-PROP}. On the level of $K$-theory 
each vertex of this diagram represents a map $A(B) \to K(\Q)$ and the edges represent 
homotopies of such maps. The outer rectangle coincides with diagram (\ref{HOMOTINV_MAPS}).  

In order to show that the diagram (\ref{HOMOTINV_MAPS}) can be filled by a homotopy of 
homotopies it is enough to show that each of the subdiagrams (1)-(3) in the above diagram 
of functors can be filled by a homotopy of homotopies. In the case of  subdiagram (1) such a homotopy
 of homotopies exists since this subdiagram commutes. By Proposition  \ref{MAPSFIBR-PROP}
subdiagram (2) commutes up to a natural chain  homotopy, so it again can be filled by 
a homotopy of homotopies.  Proposition \ref{MAPSFIBR-PROP} says also that the maps 
$f^{\infty}_{X}$ preserve the homological  filtration of the twisted tensor products and that they 
induce the map $\id\otimes (f|_{F_{1}})_{\ast}$ on the filtration  quotients. 
This, together with the fact that the map 
 $(f|_{F_{1}})_{\ast}\colon H_{\ast}(F_{1})\to H_{\ast}(F_{2})$ is an isomorphism of $\pi_{1}B$-modules, 
 implies that we also have a homotopy of homotopies filling subdiagram (3).

\end{proof}

%%%%%%%%%%%%%%%%%%%%%%%%%%%%
% ADDITIVITY OF THE SECONDARY TRANSFER
%%%%%%%%%%%%%%%%%%%%%%%%%%%%

\section{Additivity of the secondary transfer}
\label{ADDITIVITY SEC}

Our goal of this section is to prove that  secondary transfer maps have additivity properties that are
analogous to additivity of the Becker-Gottlieb transfer and the $A$-theory transfer. 
We start by considering additivity of the homotopy secondary transfer:

\begin{theorem}
\label{HADD THM}
For $i=0, 1, 2$ let $p_{i}\colon E_{i}\to B$ be a unipotent fibration. Assume that we have maps 
of fibrations  
\begin{equation*}
\begin{tikzpicture}[baseline=(current bounding box.center)]
\matrix (m) 
[matrix of math nodes, row sep=3em, column sep=2.5em, text height=1.5ex, text depth=0.25ex]
{
E_{1} & E_{0}&  E_{2} \\
& B &  \\
};
\path[->, thick, font=\scriptsize]
(m-1-1)
edge node[anchor= north east] {$p_{1}$} (m-2-2)
(m-1-2)
edge node[anchor=  east] {$p_{0}$} (m-2-2)
edge node[anchor= south] {} (m-1-1)
edge node[anchor= south] {j} (m-1-3)
(m-1-3)
edge node[anchor= north west] {$p_{2}$}  (m-2-2)
; 
\end{tikzpicture}
\end{equation*}
where $j$ is a cofibration over $B$. Let $E:= E_{1}\cup_{E_{0}}E_{2}$ and let 
$p\colon E\to B$ be the fibration given by $p:=p_{1}\cup_{p_{0}}p_{2}$. Then $p$ is 
a unipotent fibration and we have 
\begin{equation}
\label{HADD EQ}
[\wh(p^{!})] = [k_{1\ast}\wh(p_{1})] + [ k_{2\ast}\wh(p_{2})] - [k_{0\ast}\wh(p_{0})] 
\end{equation}
Here $k_{i\ast}\colon \wh(E_{i})\to \wh(E)$ is induced by the map $k_{i}\colon E_{i}\to E$. 
\end{theorem}

\begin{lemma}
\label{HOMOTCOMM LEMMA}
Consider a diagram of chain complexes
\begin{equation*}
\begin{tikzpicture}[baseline=(current bounding box.center)]
\matrix (m) 
[matrix of math nodes, row sep=2em, column sep=2em, text height=1.5ex, text depth=0.25ex]
{
A & B \\
A' &  B'  \\
};
\path[->, thick, font=\scriptsize]
(m-1-1)
edge node[auto] {$f$} (m-1-2)
edge node[anchor=east] {$g$} (m-2-1)
(m-1-2)
edge node[anchor=west] {$g'$}  (m-2-2)
(m-2-1) 
edge  node[anchor=south] {$f'$} (m-2-2)
; 
\end{tikzpicture} 
\end{equation*}
that commutes up to a chain homotopy $h$.
There exists a map $\tilde g'\colon \Cyl(f)\to B'$
such that the diagram 
\begin{equation*}
\begin{tikzpicture}[baseline=(current bounding box.center)]
\matrix (m) 
[matrix of math nodes, row sep=2em, column sep=2em, text height=1.5ex, text depth=0.25ex]
{
A & \Cyl(f) \\
A' &  B'  \\
};
\path[->, thick, font=\scriptsize]
(m-1-1)
edge node[auto] {} (m-1-2)
edge node[anchor=east] {$g$} (m-2-1)
(m-1-2)
edge node[anchor=west] {$\tilde g'$}  (m-2-2)
(m-2-1) 
edge  node[anchor=south] {$f'$} (m-2-2)
; 
\end{tikzpicture} 
\end{equation*}
commutes.  Moreover $\tilde g'$ is chain homotopic to the composition 
$$\Cyl(f) \to B \overset{g'}{\lra} B'$$
\end{lemma}

\begin{proof}
Recall that $\Cyl(f)_{n} = A_{n}\oplus A_{n-1}\oplus B_{n}$. 
The map $\tilde g'_{n}\colon \Cyl(f)_{n}\to B'_{n}$ is given by 
$$\tilde g'_{n}(a_{1}, a_{2}, b) = f'g(a_{1}) + h(a_{2}) + g'(b)$$
The second statement of the lemma is easy to verify. 
\end{proof}

\ 

\begin{proof}[Proof of Theorem \ref{HADD THM}]
Let $F$ denote the fiber of the fibration $p$ and for  $i=0, 1, 2$ let $F_{i}$ be the fiber of $p_{i}$. 
The fact that the action of $\pi_{1}(B)$ on  $F= F_{1}\cup_{F_{0}}F_{2}$ is unipotent can be obtained 
using the Mayer-Vietoris sequence for homology of the fibers. 

The strategy of the proof of additivity of the secondary transfer is as follows. We will construct maps 
 $f^{h}_{1}\colon \wh(B)\to \wh(E)$ and  $f^{h}_{2}\colon \wh(B)\to \wh(E_{2})$ such that 
$[f^{h}_{1}] = [k_{2\ast}f^{h}_{2}] $. We will also show that 
$$[\wh(p^{!})] = [k_{1\ast}\wh(p_{1}^{!})] + [f_{1}^{h}] \ \ \ \ \ \text{and} \ \ \ \ \ 
[\wh(p_{2}^{!})] = [j_{\ast}\wh(p_{0}^{!})] + [f_{2}^{h}]$$ 
Since $k_{0} = k_{2}j$ the second of these equations will give 
$$[k_{2\ast}\wh(p_{2}^{!})] = [k_{0\ast}\wh(p_{0}^{!})] + [k_{2\ast} f_{2}^{h}] =  [k_{0\ast}\wh(p_{0}^{!})] + [f_{1}^{h}] $$
which, combined with the first equation, will yield the formula (\ref{HADD EQ}). 

The construction of  the map $f^{h}_{1}$ will proceed following the  scheme outlined in 
(\ref{HOMOTFIBMAP NN}). First, we will construct a map $f_{1}^{A}\colon A(B)\to A(E)$. 
Subsequently we will consider the diagram 
\begin{equation}
\label{H1 EQ}
\begin{tikzpicture}[baseline=(current bounding box.center)]
\matrix (m) 
[matrix of math nodes, row sep=3em, column sep=5em, text height=1.5ex, text depth=0.25ex]
{
A(B) & A(E) \\
K(\Q) & K(\Q) \\ 
};
\path[->, thick, font=\scriptsize]
(m-1-1) 
edge node[auto] {$f_{1}^{A}$} (m-1-2)
edge node[anchor=east] {$\lambda^{h}_{B}$} (m-2-1)
(m-2-1) 
edge node[below] {$\chi(\Cone(k_{1}|_{F_{1}\ast}))$} (m-2-2)
(m-1-2)
edge node[auto] {$\lambda^{h}_{E}$} (m-2-2)
; 
\end{tikzpicture}
\end{equation}
Here $\Cone(k_{1}|_{F_{1}\ast})$ denotes the mapping cone of the map 
$k_{1}|_{F_{1}\ast}\colon H_{\ast}(F_{1}) \to H_{\ast}(F)$, and the map
$$\chi(\Cone(k_{1}|_{F_{1}\ast}))\colon K(\Q)\to K(\Q)$$ 
is  induced by the  functor $\Chfd(\Q)\to \Chfd(\Q)$ given by tensoring by 
$\Cone(k_{1}|_{F_{1}\ast})$. We will show that the diagram (\ref{H1 EQ}) commutes up to 
a preferred homotopy $h_{1}$. This homotopy together with the map $f_{1}^{A}$ will 
define the map $f^{h}_{1}$. 

In order to obtain the map $f_{1}^{A}$ recall that the map $k_{1\ast}A(p_{1}^{!})$ is induced by 
the functor that assigns to a space $X\in\Rfd(B)$ the space 
$k_{1\ast}p_{1}^{\ast}X\in \Rfd(E)$, and that the map $A(p^{!})$
is induced by the functor $X\mapsto p^{\ast}X$. For $X\in \Rfd(B)$  we have a cofibration 
$k_{1\ast}p_{1}^{\ast}X \to p^{\ast}X$. Let $M_{X} \in \Rfd(E)$ 
denote the cofiber of this map. The assignment $X\mapsto M_{X}$ defines an exact functor
$\Rfd(B)\to \Rfd(E)$. The map $f_{1}^{A}$ is induced by this functor.   
The above constructions give a short exact sequence of functors 
\begin{equation}
\label{MXSES EQ}
k_{1\ast}p_{1}^{\ast}X\to p^{\ast}X \to M_{X}
\end{equation}
Applying Waldhausen's additivity theorem we obtain a homotopy 
\begin{equation}
\label{FA1 EQ}
A(p^{!}) \simeq k_{1\ast}A(p_{1}^{!}) + f_{1}^{A}
\end{equation}
Next, in order to describe a homotopy that fills the diagram (\ref{H1 EQ}) consider the following 
diagram of functors:
\begin{equation}
\label{ADDFUNCT EQ}
\begin{tikzpicture}[baseline=(current bounding box.center)]
\matrix (m) 
[matrix of math nodes, row sep=3em, column sep=1em, text height=1.5ex, text depth=0.25ex, font=\small]
{
C_{\ast}(k_{1\ast}p_{1}^{\ast}X, E) & & C_{\ast}(M_{X}, E) \\
C_{\ast}(p_{1}^{\ast}X, E_{1}) & C_{\ast}(p^{\ast}X, E) & \coker(k_{X\ast})\\
C_{\ast}(p_{1}^{\ast}X, E_{1}) & \Cyl(k_{X\ast}) & \Cone(k_{X\ast})\\
\bC(X, B)\otimes_{p_{1,X}}H_{\ast}(F_{1}) & \Cyl(k_{X}^{\infty}) & \Cone(k_{X}^{\infty})\\
C_{\ast}(X, B)\otimes H_{\ast}(F_{1}) & C_{\ast}(X, B)\otimes\Cyl(k_{1}|_{F_{1}\ast})
&  C_{\ast}(X, B)\otimes\Cone(k_{1}|_{F_{1}\ast})\\
};

\path[->, thick, font=\scriptsize]
(m-1-1) 
edge node[auto] {} (m-2-2)
(m-2-1)
edge node[auto] {$\simeq$} (m-1-1)
edge node[anchor=north ] {$k_{X\ast}$} (m-2-2)
(m-2-2)
edge  (m-1-3)
edge node[anchor=north ] {} (m-2-3)
(m-2-3)
edge node[anchor=west] {$\simeq$}(m-1-3)
(m-3-1)
edge  node[anchor=east] {$=$} (m-2-1)
edge node[anchor= south] {} (m-3-2)
edge  node[anchor= east] {$\beta_{p_{1,X}}$} (m-4-1)
(m-3-2)
edge  (m-2-2)
edge node[anchor= south] {} (m-3-3)
edge  node[anchor= west] {$g_{X}$} (m-4-2)
(m-3-3)
edge  (m-2-3)
edge  (m-4-3)
(m-4-1)
edge node[anchor= south] {} (m-4-2)
(m-4-2)
edge node[anchor= south] {} (m-4-3)
(m-5-1)
edge  (m-5-2)
(m-5-2)
edge  (m-5-3)
;
\path[solid,  thick, font=\scriptsize]  
(m-4-1) 
edge node[anchor=east] {additivity} (m-5-1)
(m-4-2) 
edge node[anchor=west] {additivity} (m-5-2)
(m-4-3) 
edge node[anchor=west] {additivity} (m-5-3)
; 
\end{tikzpicture}
\end{equation}
The map $k_{X\ast}$ in this diagram is induced by the map of fibrations 
$k_{X}\colon p_{1}^{\ast}X\to p^{\ast}X$. The complexes $\Cyl(k_{X\ast})$
and $\Cone(k_{X\ast})$ are, respectively, the mapping cylinder and the mapping 
cone of $k_{X\ast}$. Similarly $\Cyl(k^{\infty}_{X})$ and $\Cone(k^{\infty}_{X})$ are 
the mapping cylinder and the mapping cone of the map 
$$k_{X}^{\infty}\colon \bC(X, B)\otimes_{p_{1,X}} H_{\ast}(F_{1}) 
\to \bC(X, B)\otimes_{p_{X}} H_{\ast}(F)$$
given by Proposition \ref{MAPSFIBR-PROP}.  Finally, $\Cyl(k_{1}|_{F_{1}\ast})$ and 
$\Cone(k_{1}|_{F_{1}\ast})$ are the mapping cylinder and the mapping cone of the map
$k_{1}|_{F_{1}\ast}\colon H_{\ast}(F_{1})\to H_{\ast}(F)$. The horizontal maps are defined 
in the obvious way so that each row of the diagram forms a short exact sequence. 

All vertical maps are quasi-isomorphism. They are defined in the obvious way with the exception 
the map $g_{X}$ which is given as follows. By Proposition \ref{MAPSFIBR-PROP} we have 
a diagram
\begin{equation*}
\begin{tikzpicture}[baseline=(current bounding box.center)]
\matrix (m) 
[matrix of math nodes, row sep=4em, column sep=2em, text height=1.5ex, text depth=0.25ex]
{
C_{\ast}(p_{1}^{\ast}X, E_{1}) & C_{\ast}(p^{\ast}X, E) \\
\bC(X, B)\otimes_{p_{1X}} H_{\ast}(F_{1}) &   \bC(X, B)\otimes_{p_{X}} H_{\ast}(F)  \\
};
\path[->, thick, font=\scriptsize]
(m-1-1)
edge node[auto] {$k_{X\ast}$} (m-1-2)
edge node[anchor=east] {$\beta_{p_{1, X}}$} (m-2-1)
(m-1-2)
edge node[anchor=west] {$\beta_{p_{X}}$}  (m-2-2)
(m-2-1) 
edge  node[anchor=south] {$k_{X}^{\infty}$} (m-2-2)
; 
\end{tikzpicture} 
\end{equation*}
that commutes up to a chain homotopy. Using Lemma  \ref{HOMOTCOMM LEMMA} we obtain 
from here a commutative diagram 
\begin{equation*}
\begin{tikzpicture}[baseline=(current bounding box.center)]
\matrix (m) 
[matrix of math nodes, row sep=4em, column sep=2em, text height=1.5ex, text depth=0.25ex]
{
C_{\ast}(p_{1}^{\ast}X, E_{1}) & \Cyl(k_{X\ast}) \\
\bC(X, B)\otimes_{p_{1X}} H_{\ast}(F_{1}) &   \bC(X, B)\otimes_{p_{X}} H_{\ast}(F)  \\
};
\path[->, thick, font=\scriptsize]
(m-1-1)
edge node[auto] {} (m-1-2)
edge node[anchor=east] {$\beta_{p_{1, X}}$} (m-2-1)
(m-1-2)
edge node[anchor=west] {$\bar{\beta}_{p_{X}}$}  (m-2-2)
(m-2-1) 
edge  node[anchor=south] {$k_{X}^{\infty}$} (m-2-2)
; 
\end{tikzpicture} 
\end{equation*}
The map $g_{X}$ is the composition of $\bar{\beta}_{p_{X}}$ and the inclusion 
$$\bC(X, B)\otimes_{p_{X}} H_{\ast}(F) \to \Cyl(k_{X}^{\infty})$$
The lowest vertical edges in  the diagram indicate a passage between chain complexes 
using additivity. This means the following construction. Since the map $k^{\infty}_{X}$
preserves the homological filtrations on  $ \bC(X, B)\otimes_{p_{1 X}} H_{\ast}(F_{1})$
and  $\bC(X, B)\otimes_{p_{X}} H_{\ast}(F)$ the complexes $\Cyl(k^{\infty}_{X})$
and $\Cone(k^{\infty}_{X})$ are endowed with induced filtrations. Each lowest vertical 
edge indicates a passage from the filtered chain complex to the direct sum of the filtration 
quotients. As usual, after passage to the induced maps $A(B)\to K(\Q)$ each additivity edge 
gives a homotopy obtained using Waldhausen's additivity theorem. We note here that the 
additivity edges are related to one another as follows. The maps in the short exact sequence 
$$\bC(X, B)\otimes_{p_{1, X}}H_{\ast}(F_{1})\lra \Cyl(k^{\infty}_{X}) \lra \Cone(k^{\infty}_{X})$$
preserve filtrations. Moreover, their restrictions to the filtration subcomplexes also form  
short exact sequences, and so do the induced maps of filtrations quotients. The bottom 
row of the diagram is the direct sum of the short exact sequences of these filtration quotients. 

Each vertex of in the diagram (\ref{ADDFUNCT EQ}) induces a map $A(B)\to K(\Q)$. All vertical edges
define homotopies of such maps. Concatenation of homotopies defined by rightmost vertical edges 
gives a homotopy filling the diagram (\ref{H1 EQ}).  This homotopy, together with the map $f_{1}^{A}$
defines the map $f^{h}_{1}\colon \wh(B)\to \wh(E)$. 

The existence of a homotopy between the maps $\wh(p^{!})$ and $k_{1\ast}\wh(p_{1}^{!}) + f^{h}_{1}$
follows directly from the above construction. The map $k_{1\ast}\wh(p_{1}^{!})$ is defined by the map 
$k_{1\ast}A(p_{1}^{!})$ and the homotopy induced by the leftmost vertical edges in the diagram 
(\ref{ADDFUNCT EQ}). The map $\wh(p^{!})$ is homotopic to the map defined by $A(p^{!})$
and the homotopy induced by the middle vertical edges in (\ref{ADDFUNCT EQ}). As we have already 
noticed applying Waldhausen's additivity theorem to the short exact sequence of functors 
(\ref{MXSES EQ}) defines a homotopy between $A(p^{!})$ and $k_{1\ast}A(p_{1}^{!}) + f_{1}^{A}$. 
In order to lift this homotopy to a homotopy between $\wh(p^{!})$ and $k_{1\ast}\wh(p_{1}^{!}) + f^{h}_{1}$
it is enough to apply the additivity theorem to the horizontal short exact sequences in the diagram 
(\ref{ADDFUNCT EQ}).

Construction of the map $f^{h}_{2}\colon \wh(B)\to \wh(E_{2})$ proceeds in  the same way as 
the construction of $f^{h}_{1}$, with the cofibration $j\colon E_{0}\to E_{2}$ used in place of $k_{1}$. 
By the same argument as above we obtain a homotopy $\wh(p_{2}^{!}) \simeq j_{\ast}\wh(p_{0}^{!}) + f^{h}_{2}$.
Finally, the fact that the maps $f^{h}_{1}$ and $k_{2\ast}f^{h}_{2}$ are homotopic can be verified directly 
by inspecting  the construction of $f^{h}_{1}$ and $f^{h}_{2}$.
\end{proof}

A statement analogous to Theorem \ref{HADD THM} holds for the smooth secondary transfer.  
Given a smooth bundle $p\colon E\to B$ whose fibers are manifolds with a boundary by 
the vertical boundary of $p$ we will understand the smooth bundle $\partial^{v}p \colon \partial^{v}E \to B$
obtained by restricting $p$ to the union of boundaries of its fibers.  We have: 

\begin{theorem}
\label{SADD THM}
Let $p\colon E\to B$ be a smooth  bundle with closed fibers, and for $i=0, 1, 2$ let  
$p_{i}\colon E_{i}\to B$ be unipotent subbundles of $p$ such that $p_{0}$ is the vertical boundary 
of both $p_{1}$ and $p_{2}$, and that $E = E_{1}\cup_{E_{0}} E_{2}$.  Then $p$ is a unipotent bundle
and we have 
\begin{equation*}
\label{SADD EQ}
[\wh[s](p^{!})] = [k_{1\ast}\wh[s](p_{1})] + [ k_{2\ast}\wh[s](p_{2})] - [k_{0\ast}\wh[s](p_{0})] 
\end{equation*}
Here $k_{i\ast}\colon \wh[s](E_{i})\to \wh[s](E)$ is induced by the map $k_{i}\colon E_{i}\to E$. 
\end{theorem}

\begin{note}
Recall that the construction of the smooth secondary transfer we are working with 
uses the combinatorial model of $Q(E_{+})$  built using partitions (\ref{AKQ CONSTR NN}).  
This model is functorial with respect 
to embeddings of submanifolds of codimension $0$. As a consequence in the notation 
of Theorem \ref{SADD THM} for $i=1, 2$ the inclusion maps $k_{i}\colon E_{i} \to E$ induce maps 
$k_{i\ast}\colon Q(E_{i+}) \to Q(E)$, which then lift to maps  $k_{i\ast}\colon \wh[s](E_{i})\to \wh[s](E)$.  
The map $k_{0\ast}\colon Q(E_{0+}) \to Q(E_{+})$ is constructed as follows. 
Let  $b\colon E_{0}\times [-1, 1] \to E$ be a fiberwise bicollar neighborhood of $E_{0}$. Thus, 
$b$ is a smooth embedding of bundles over $B$ such that $b(E_{0}\times \{0\}) = E_{0}$, 
$b(E_{0}\times [-1, 0]) \subseteq E_{1}$ and $b(E_{0}\times [0, 1]) \subseteq E_{2}$. 
The inclusion $k_{0}\colon E_{0}\to E$ coincides with the composition 
$$E_{0}  \to E_{0}\times \{0\} \subseteq E_{0}\times [-1, 1] \overset{b}{\lra} E$$
For a partition $P\subseteq E_{0}\times I$ the submanifold $P\times [-1, 1] \subseteq E_{0}\times [-1, 1]$ 
defines (modulo permutation of coordinates) a partition of $(E_{0}\times [-1, 1]) \times I$. 
This assignment induces a map $Q(E_{0+}) \to Q(E_{0}\times [-1, 1]_{+})$
Furthermore, since $b$ is an embedding of codimension $0$ it gives a map 
$b_{\ast}\colon  Q(E_{0}\times [-1, 1]_{+}) \to Q(E_{+})$.  Composing these two maps we obtain 
a map $k_{0\ast} \colon Q(E_{0+}) \to Q(E_{+})$ which lifts to 
$$k_{0\ast}\colon \wh[s](E_{0})\to \wh[s](E)$$ 
\end{note}

\begin{proof}[Proof of Theorem \ref{SADD THM}]
The basic scheme of the proof is the same as that of the proof of Theorem \ref{HADD THM}:
it suffices to show that there exist maps $f^{s}_{1}\colon \wh[s](B)\to \wh[s](E)$ and  
$f^{s}_{2}\colon \wh[s](B)\to \wh[s](E_{2})$ such that 
$$[\wh[s](p^{!})] = [k_{1\ast}\wh[s](p_{1}^{!})] + [f^{s}_{1}] \ \ \ \ \ \ \ \ \ \ 
[\wh[s](p_{2}^{!})] = [j_{\ast}\wh[s](p_{0}^{!})] + [f^{s}_{2}]$$ 
where $j\colon E_{0} \to E_{2}$ is the inclusion map, 
and that $[f^{s}_{1}] = [k_{2\ast}f^{s}_{2}]$. 
Moreover, the construction of the maps $f^{s}_{1}$, $f_{2}^{s}$ and verification that they satisfy the 
above identities also mimics  the arguments we used in the proof of Theorem \ref{HADD THM}.  
Recall  that the map $f^{h}_{1}$ in that proof was defined using a map $f^{A}_{1}\colon A(B) \to A(E)$
and a preferred homotopy in the diagram (\ref{H1 EQ}). Similarly, in order to define the map $f^{s}_{1}$
we need to construct a map $f^{Q}_{1}\colon Q(B_{+}) \to Q(E_{+})$ together 
with a preferred homotopy $\lambda_{E} f_{1}^{Q} \simeq \chi(\Cone(k_{1}|F_{1\ast}))\lambda_{B}$. 
Furthermore, in order to obtain a homotopy 
$\wh[s](p^{!}) \simeq  k_{1\ast}\wh[s](p_{1}^{!}) + f^{s}_{1}$
it suffices to construct a homotopy 
\begin{equation}
\label{F1QHOMOTOPY EQ}
Q(p^{!}) \simeq k_{1\ast}Q(p_{1}^{!}) +f_{1}^{Q}
\end{equation}
together with an appropriate homotopy of homotopies. 

The construction of the map $f_{1}^{Q}$ can be simplified in two ways. First, let $\SB$  denote the simplicial 
set of smooth singular simplices of $B$. It will suffice to construct a map $f_{1}^{Q}\colon |\SB| \to Q(E_{+})$. 
Since   $B\simeq |\SB|$ this map will extend to a map of infinite loop spaces $f_{1}^{Q}\colon Q(B_{+}) \to Q(E_{+})$. 
Second, the construction of the map $f_{1}^{Q}\colon |\SB| \to Q(E_{+})$
can be reduced to the following combinatorial construction. Let $\mathcal{T}_{1}\mathcal{P}_{\bullet}(E\times I)$
denote the simplicial category given by the first stage of the $\mathcal{T}_{\bullet}$-construction. Thus for $k\geq 0$
the objects of $\mathcal{T}_{1}\mathcal{P}_{k}(E\times I)$ are pairs $(P_{0}, P_{1})$ where $P_{1}$
is a  bundle of partitions of $E\times I$ over the standard simplex $\Delta^{k}$, and $P_{0}$
is a subbundle of $P_{1}$. Considering $\SB$ as a (discrete) simplicial category it will suffice to give a functor 
$F_{1}^{Q} \colon \SB \to \mathcal{T}_{1}\mathcal{P}_{\bullet}(E\times I)$. Such functor  will determine the map 
$f_{1}^{Q} \colon |\SB| \to Q(E_{+})$ (cf. \cite[2.5]{BDW}). 

In order to describe the functor $F_{1}^{Q}$ we will use the setup of \cite[4.3]{BDW}. Using the notation 
introduced there for a  $k$-simplex $\sigma\in \SB$  by $\sigma^{\ast}TB_{\varepsilon}$  we denote the disc bundle over 
$\Delta^{k}$ induced by $\sigma$ from the disc bundle of the tangent bundle of $B$. 
The exponential map defines a map of bundles 
\begin{equation*}
\begin{tikzpicture}[baseline=(current bounding box.center)]
\matrix (m) 
[matrix of math nodes, row sep=2em, column sep=1.5em, text height=1.5ex, text depth=0.25ex]
{
\sigma^{\ast}TB_{\varepsilon} & &  B\times \Delta^{k} \\
& \Delta^{k} &  \\
};
\path[->, thick, font=\scriptsize]
(m-1-1)
edge  (m-2-2)
edge node[anchor= south] {$\exp$} (m-1-3)
(m-1-3)
edge  (m-2-2)
; 
\end{tikzpicture}
\end{equation*}
which is a fiberwise embedding. 
Let $\frac{1}{3} < a < b < 1$. Define
$$P_{\sigma} := \sigma^{\ast}TB_{\varepsilon}\times [a, b] \cup (B\times  \Delta^{k}) \times [0, \tfrac{1}{3}]$$ 
This space admits a map  $P_{\sigma} \to B\times  \Delta^{k}\times I$ which lets us   consider it
(modulo permutation of factors) as a bundle of partitions of $B\times I$ over $\Delta^{k}$. 
For $\sigma \in \SB$ we set
$$F_{1}^{Q}(\sigma) := ((p_{1}\times \id)^{-1}(P_{\sigma}) \cup E\times \Delta^{k}\times [0, \tfrac{1}{3}], \ 
(p\times \id)^{-1}(P_{\sigma}))$$

In order to describe the homotopy (\ref{F1QHOMOTOPY EQ})
we will use Waldhausen's pre-additivity theorem \cite[3.6]{BDKW}. Take 
$\mathcal{T}_{2}\mathcal{P}_{\bullet}(E\times I)$ to be the second stage of the $\mathcal{T}_{\bullet}$-construction 
on $\mathcal{P}_{\bullet}(E\times I)$. It is a simplicial category with objects in the $k$-th simplicial dimension
given by triples $(P_{0}, P_{1}, P_{2})$  where $P_{i}$ is a bundle of partitions of $E\times I$
over $\Delta^{k}$, and  $P_{i}$ is a subbundle of $P_{i+1}$. For $0\leq i < j \leq j$ let 
$D_{ij}\colon \mathcal{T}_{2}\mathcal{P}_{\bullet}(E\times I) \to \mathcal{T}_{1}\mathcal{P}_{\bullet}(E\times I)$
denote the functor given by $D_{ij}(P_{0}, P_{1}, P_{2}) = (P_{i}, P_{j})$. For a functor 
$H\colon \SB \to \mathcal{T}_{2}\mathcal{P}_{\bullet}(E\times I)$ the preadditivity theorem gives a homotopy 
between the map $|\SB| \to Q(E_{+})$ induced by the functor $D_{02}H$ and the sum of maps induced by
the functors $D_{01}H$ and $D_{12}H$.  Consider the functor $H$ given by 
$$H(\sigma) := (E\times \Delta^{k}\times [0, \tfrac{1}{3}],  \ 
(p_{1}\times \id)^{-1}(P_{\sigma}) \cup E\times \Delta^{k}\times [0, \tfrac{1}{3}], \ 
(p\times \id)^{-1}(P_{\sigma}))$$
Notice that $D_{12}H = F^{Q}_{1}$. It remains to notice that by \cite[\S 4]{BDW} the maps 
induced by the functors $D_{01}H$ and $D_{02}H$ are, respectively, $k_{1\ast}Q(p_{1}^{!})\eta$ 
and $Q(p^{!})\eta$ where $\eta \colon |\SB| \to Q(B_{+})$ is the coaugmentation map.  

The above constructions of the map $f_{1}^{Q}$ and the homotopy $Q(p^{!}) \simeq Q(p_{1}^{!}) + f_{1}^{Q}$ 
are  obtained by replicating the constructions of the maps $f_{1}^{h}$ and the homotopy 
$A(p^{!}) \simeq A(p_{1}^{!}) + f_{1}^{A}$ from the proof of Theorem \ref{HADD THM}, using partitions and the 
$\mathcal{T}_{\bullet}$-construction in place of retractive spaces and the $\mathcal{S}_{\bullet}$-construction.
The arguments that show that the map $f_{1}^{Q}$ admits a lift to $f_{1}^{s}\colon \wh[s](B) \to \wh[s](E)$
and that the homotopy (\ref{F1QHOMOTOPY EQ}) lifts to a homotopy 
$\wh[s](p^{!}) \simeq  k_{1\ast}\wh[s](p_{1}^{!}) +  f^{s}_{1}$ involve constructions on the level of chain 
complexes that duplicate the constructions from the proof of Theorem \ref{HADD THM}. 
The map $f_{2}^{s}$ and the homotopy $\wh[s](p_{2}^{!}) = j_{\ast}\wh[s](p_{0}^{!}) + f^{s}_{2}$ are obtained 
in an analogous way. 

\end{proof}

%%%%%%%%%%%%%%%%%%%%%%%%%%%%
% COMPOSITION OF UNIPOTENT FIBRATIONS
%%%%%%%%%%%%%%%%%%%%%%%%%%%%

\section{Composition of unipotent fibrations}
\label{COMPUNIP SEC}

Our next goal is to give a proof of Theorem \ref{MAIN2 THM}. Recall that this theorem states 
that the homotopy secondary transfer preserves compositions of unipotent fibrations, and
the smooth secondary transfer preserves compositions of unipotent bundles. 

The statement of Theorems \ref{MAIN2 THM} relies 
on the fact that the composition of unipotent bundles (or unipotent fibrations) 
is again unipotent. While this property 
is implicitly present in the work of Igusa \cite{IgusaAx} we give its proof below for completeness, 
and also because its main ingredient,  Lemma \ref{IGUSA-LEMMA},  will be needed 
later on.

\

\begin{lemma}[{\cite[Lemma 8.9]{IgusaAx}}]
\label{IGUSA-LEMMA}
Let 
$$F_{p}\to E \overset{p}{\lra} B$$ 
be a unipotent fibration. 
There is a finite sequence of unipotent fibrations
\begin{equation*}
\begin{tikzpicture}[baseline=(current bounding box.center)]
\matrix (m) 
[matrix of math nodes, row sep=3em, column sep=2.5em, text height=1.5ex, text depth=0.25ex]
{
E_{0}& E_{1} & \dots & E_{k} \\
& B &   &  \\
};
\path[->, thick, font=\scriptsize]
(m-1-1)
edge  (m-1-2)
edge node[anchor=north east] {$p_{0}$} (m-2-2)
(m-1-2)
edge  (m-1-3)
edge node[anchor= west] {$p_{1}$} (m-2-2)
(m-1-3)
edge  (m-1-4)
(m-1-4)
edge node[anchor=north west] {$p_{k}$} (m-2-2.north east)
; 
\end{tikzpicture}
\end{equation*}
such that 
\begin{enumerate}
\item[(i)]  $E_{0} = \Sigma^{n}_{B}E$ for some $n\geq 0$;
\item[(ii)] $p_{k}\colon E_{k} \to B$ is a rational homotopy equivalence;
\item[(iii)] for every $i$ we have a cofibration sequence over $B$:
$$B\times S^{n_{i}} \overset{\alpha_{i}}{\lra} E_{i} \lra E_{i+1}$$
(i.e.  $E_{i+1} = E_{i}\cup_{B\times S^{n_{i}}}B\times D^{n_{i}+1}$).
\end{enumerate}
\end{lemma}

We denote here by  $\Sigma_{B}^{n}E \to B$  the $n$-fold fiberwise suspension of the fibration 
$p$, while $B\times S^{n_{i}}\to B$ and $B\times D^{n_{i}+1}\to B$ are, respectively, a product sphere bundle
and a product disc bundle.    

\ 

\begin{theorem}[Igusa] If $q\colon D\to E$ and $p\colon E\to B$ are unipotent bundles
(resp. unipotent fibrations) then the composition $pq\colon D\to B$ is also a unipotent 
bundle (resp. a unipotent fibration). 
\end{theorem}

\begin{proof}
Let $F_{p}$, $F_{q}$, $F_{pq}$ denote the fibers of $p$, $q$, $pq$, respectively.
Since the only non-trivial property we need to verify is that the action of $\pi_{1}B$ on 
$H_{\ast}(F_{pq})$ is unipotent it is enough to show that the statement of the  
theorem  holds for unipotent fibrations. 
We will split our argument into a few steps. 

\emph{Step 1.} The fibration $pq$ is unipotent for an arbitrary unipotent 
fibration $p$ and any product fibration $q\colon E\times F_{q}\to E$.  

Indeed, in this case we have an isomorphism of $\pi_{1}B$-modules 
$$H_{\ast}(F_{pq}) \cong H_{\ast}(F_{p})\otimes H_{\ast}(F_{q})$$
where the action of $\pi_{1}B$ on the right hand side is given by
$\alpha(x\otimes y) = \alpha x\otimes y$. 

\emph{Step 2.} The fibration $pq$ is unipotent for an arbitrary unipotent 
fibration $p\colon E\to B$ and any fibration $q\colon D\to E$ with fiber $F_{q}$ 
such that  $\tilde H_{\ast}(F_{q}) = 0$. 

This holds since the map $(q|_{F_{pq}})_{\ast} \colon H_{\ast}(F_{pq})\to H_{\ast}(F_{p})$ is 
in this case an isomorphism of $\pi_{1}B$-modules.

\emph{Step 3.} For  $i=0, 1, 2$ let $p_{i}\colon E_{i}\to B$ be a fibration with 
a path connected base space  $B$ and 
fiber $F_{i}$ of a finite homotopy type.
 Assume that we have maps of fibrations
$$E_{1} \overset{f}{\lla} E_{0} \overset{i}{\lra} E_{2}$$
where $i$ is a cofibration over $B$. Let 
$p \colon E_{1}\cup_{E_{0}}E_{2} \to B$
be the pushout map. If three of the fibrations $p_{0}, p_{1}, p_{2}, p$ are unipotent 
then so is the fourth. 

This follows from the Mayer-Vietoris sequence for the homology of the fibers and the fact that unipotent 
$\pi_{1}B$-modules form a Serre category. 

\emph{Step 4.} As an application of Step 3 we obtain that 
if $p, \ q$ are unipotent fibrations and $\Sigma_{E}q\colon \Sigma_{E}D \to E$
is a fiberwise suspension of $q$ then $pq$ is a unipotent fibration if and only if $p\Sigma_{E}q$ is unipotent.

\emph{Step 5.} Assume now that $p$ and $q$ are arbitrary unipotent fibrations. Applying 
Lemma \ref{IGUSA-LEMMA} to $q$ we obtain a sequence of fibrations 
\begin{equation*}
\begin{tikzpicture}[baseline=(current bounding box.center)]
\matrix (m) 
[matrix of math nodes, row sep=3em, column sep=2.5em, text height=1.5ex, text depth=0.25ex]
{
D_{0}& D_{1} & \dots & D_{k} \\
& E &   &  \\
};
\path[->, thick, font=\scriptsize]
(m-1-1)
edge  (m-1-2)
edge node[anchor=north east] {$q_{0}$} (m-2-2)
(m-1-2)
edge  (m-1-3)
edge node[anchor= west] {$q_{1}$} (m-2-2)
(m-1-3)
edge  (m-1-4)
(m-1-4)
edge node[anchor=north west] {$q_{k}$} (m-2-2.north east)
; 
\end{tikzpicture}
\end{equation*}
We will show that $pq_{i}$ is a unipotent fibration for all $i$. In case 
of $i= k$ this is a consequence of Step 2.  Arguing inductively, assume that $pq_{i+1}$ is 
unipotent for some $i$. Since $pq_{i+1}$ is the pushout in the diagram of fibrations over B
$$E\times D^{n_{i}+1} \la E\times S^{n_{i}} \ra D_{i}$$
and the fibrations $E\times D^{n_{i}+1}\to B$, $E\times S^{n_{i}}\to B$ are unipotent by Step 1,
we obtain using Step 3 that $pq_{i}$ is also a unipotent fibration. As a consequence 
we get that $pq_{0}$ is  unipotent. Since $q_{0}$ is an iterated fiberwise 
suspension of $q$ by Step 4 we obtain that $pq$ is unipotent as well. 
\end{proof}

%%%%%%%%%%%%%%%%%%%%%%%%%%%%
% COMPOSITION OF SECONDARY TRANSFERS
%%%%%%%%%%%%%%%%%%%%%%%%%%%%

\section{Composition of secondary transfers}
\label{COMPOSITION  SEC}

In this section we prove Theorems \ref{MAIN2 THM} as well as its analog for 
the homotopy secondary transfer:

\begin{theorem}
\label{HMAIN2-THM}
If $p\colon E\to B$ and $q\colon D\to E$ are unipotent fibrations  then 
$$\wh[h]((pq)^{!})\simeq \wh[h](q^{!})\circ \wh[h](p^{!})$$
\end{theorem}

Our strategy will be as follows. First, in (\ref{SECCOMPRATTRIV-LEMMA})
and (\ref{SECCOMPPROD-LEMMA}) we show that 
Theorem \ref{HMAIN2-THM}  holds in two special cases, 
and then we will use an argument of Igusa to show that 
its  general statement follows from these special cases. 
Finally, we will show that essentially the same reasoning can be used to obtain 
a proof of Theorem \ref{MAIN2 THM}.

\begin{lemma}
\label{SECCOMPRATTRIV-LEMMA}
Let  $p\colon E\to B$,  $q\colon D \to E$ be unipotent fibrations  with fibers $F_{p}$
and $F_{q}$ respectively. If $\tilde H_{\ast}(F_{q}) = 0$  then 
$$\wh((pq)^{!})\simeq \wh(q^{!})\circ \wh(p^{!})$$
\end{lemma}

\begin{proof}
Let $F_{pq}$ denote the fiber of $pq$. 
Recall (\ref{ASEC TRANSFER DEF})  that the map $\wh((pq)^{!})$ is defined by the diagram 
\begin{equation*}
\begin{tikzpicture}[baseline=(current bounding box.center)]
\matrix (m) 
[matrix of math nodes, row sep=3em, column sep=4em, text height=1.5ex, text depth=0.25ex]
{
A(B) &  A(D) \\
K(\Q) &  K(\Q)  \\
};
\path[->, thick, font=\scriptsize]
(m-1-1)
edge node[auto] {$A((pq)^{!})$} (m-1-2)
edge node[anchor=east] {$\lambda^{h}_{B}$} (m-2-1)
(m-1-2)
edge node[anchor=west] {$\lambda^{h}_{D}$}  (m-2-2)
(m-2-1) 
edge  node[anchor=north] {$\chi(F_{pq})$} (m-2-2)
; 
\end{tikzpicture}
\end{equation*}
which commutes up to the homotopy $\eta^{h}_{pq}$. One the other hand, the composition 
$\wh(q^{!})\wh(p^{!})$ is induced by  the diagram 
\begin{equation*}
\begin{tikzpicture}[baseline=(current bounding box.center)]
\matrix (m) 
[matrix of math nodes, row sep=3em, column sep=3em, text height=1.5ex, text depth=0.25ex]
{
A(B) &  A(E) & A(D) \\
K(\Q) &  K(\Q)  & K(\Q) \\
};
\path[->, thick, font=\scriptsize]
(m-1-1)
edge node[auto] {$A(p^{!})$} (m-1-2)
edge node[anchor=east] {$\lambda^{h}_{B}$} (m-2-1)
(m-1-2)
edge node[auto] {$A(q^{!})$} (m-1-3)
edge node[anchor=west] {$\lambda^{h}_{E}$}  (m-2-2)
(m-1-3)
edge node[anchor=west] {$\lambda^{h}_{D}$}  (m-2-3)
(m-2-1) 
edge  node[anchor=north] {$\chi(F_{p})$} (m-2-2)
(m-2-2) 
edge  node[anchor=north] {$\chi(F_{q})$} (m-2-3)
; 
\end{tikzpicture}
\end{equation*}
The left square in this diagram commutes up to the homotopy $\eta^{h}_{p}$, and the right square 
comutes up to the homotopy $\eta^{h}_{q}$. It follows that  the outer square, defining the map 
$\wh(q^{!})\wh(p^{!})$,  commutes up to the homotopy obtained by concatenating  
$\eta_{q}\circ (A(p^{!})\times \id_{I})$ with $\chi(F_{q})\eta_{p}$.
 
By (\ref{HOMOTFIBHOMOT NN}) in order to obtain  a homotopy between $\wh(q^{!})\wh(p^{!})$ and $\wh((pq)^{!})$ 
it is enough to construct the following data:

1) a homotopy $H_{A}(p, q)\colon A(B)\times I\to A(D)$ 
between $A((pq)^{!})$ and  $A(q^{!})A(p^{!})$ ; 

2)  a homotopy $H_{K}\colon K(\Q)\times I \to K(\Q)$ between 
the maps $\chi(F_{pq})$  and  $\chi(F_{q})\chi(F_{p})$; 

3) a homotopy of homotopies that fills the following diagram:
\begin{equation}
\label{PRODHOMOT2_MAPS}
\begin{tikzpicture}[baseline=(current bounding box.center)]
\matrix (m) 
[matrix of math nodes, row sep=3em, column sep=6em, text height=1.5ex, text depth=0.25ex]
{
\lambda^{h}_{D} A(q^{!})A(p^{!}) &  \lambda^{h}_{D}A((pq)^{!}) \\
\chi(F_{q})\lambda^{h}_{E}A(p^{!}) & \\
\chi(F_{q})\chi(F_{p})\lambda^{h}_{B} &  \chi(F_{pq})\lambda^{h}_{B} \\
};
\path[solid, thick, font=\scriptsize]
(m-1-1)
edge node[auto] {$\lambda^{h}_{D}H_{A}(p,q)$} (m-1-2)
edge node[anchor=east] {$\eta^{h}_{q}\circ(A(p^{!})\times \id_{I})$} (m-2-1)
(m-1-2)
edge node[anchor=west] {$\eta^{h}_{pq}$}  (m-3-2)
(m-2-1) 
edge  node[anchor=east] {$\chi(F_{q})\eta^{h}_{p}$} (m-3-1)
(m-3-1) 
edge  node[anchor=north] {$H_{K}(\lambda^{h}_{B}\times \id_{I})$} (m-3-2)
; 
\end{tikzpicture}
\end{equation}
Each vertex of this diagram  represents a map $A(B)\to K(\Q)$ and edges represent homotopies 
of such maps.

1) \emph{Construction of $H_{A}(p, q)$.}  The map $A((pq)^{!})$ comes from the functor
$\Rfd(B)\to\Rfd(D)$ that assigns to a retractive space $X$ the space $(pq)^{\ast}X$, while the 
composition $A(q^{!})A(p^{!})$ comes from the functor that sends $X$ to $q^{\ast}p^{\ast}X$. 
The canonical isomorphisms 
$$(pq)^{\ast}X\overset{\cong}{\lra}q^{\ast}p^{\ast}X$$
define a natural transformation of functors, and so they induce a homotopy between 
 $A((pq)^{!})$ and $A(q^{!})A(p^{!})$. This is the homotopy $H_{A}(p, q)$.

2) \emph{Construction of $H_{K}$.}
Recall that the maps $\chi(F_{p}), \chi(F_{q})$ and $\chi(F_{pq})$ 
are induced by functors $\Chfd(\Q)\to \Chfd(\Q)$ than tensor a chain complex $C$
by, respectively,  $H_{\ast}(F_{p}), H_{\ast}(F_{q})$, and $H_{\ast}(F_{pq})$.
As a consequence the map $\chi(F_{q})\chi(F_{p})$ is induced by the functor that 
tensors $C\in \Chfd(\Q)$ by the chain complex $H_{\ast}(F_{p})\otimes H_{\ast}(F_{q})$. 
Since $\tilde H_{\ast}(F_{q}) = 0$ we have an isomorphism 
$$H_{\ast}(F_{pq}) \overset{(q|_{F_{q}})_{\ast}}{\lra} H_{\ast}(F_{p}) 
\overset{\cong}\lra H_{\ast}(F_{p})\otimes H_{\ast}(F_{q})$$
This induces a natural isomorphism of functors 
$$-\otimes H_{\ast}(F_{pq})\Ra -\otimes (H_{\ast}(F_{p})\otimes H_{\ast}(F_{q}))$$
which, in turn, defines the homotopy $H_{K}$.

3) \emph{Construction of the homotopy of homotopies.} In order to show that the diagram 
(\ref{PRODHOMOT2_MAPS}) can be filled by a homotopy of homotopies we will first replace 
it by the underlying diagram of functors:
\begin{equation*}
\begin{tikzpicture}[, baseline=(current bounding box.center)]
\matrix (m) 
[matrix of math nodes, row sep=4em, column sep=2em, text height=1.5ex, text depth=0.25ex, font=\small ]
{
C_{\ast}(q^{\ast}p^{\ast}X, D) &  & C_{\ast}((pq)^{\ast}X, D) \\
\bC(p^{\ast}X, E)\otimes H_{\ast}(F_{q}) & & \\[-1em]
C_{\ast}(p^{\ast}X, E)\otimes H_{\ast}(F_{q}) & C_{\ast}(p^{\ast}X, E) & \\
\bC(X, B)\otimes_{p_{X}}H_{\ast}(F_{p})\otimes H_{\ast}(F_{q})  & 
\bC(X, B)\otimes_{p_{X}}H_{\ast}(F_{p}) & 
\bC(X, B)\otimes_{(pq)_{X}}H_{\ast}(F_{pq})  \\
C_{\ast}(X, B)\otimes H_{\ast}(F_{p})\otimes H_{\ast}(F_{q}) & 
C_{\ast}(X, B)\otimes H_{\ast}(F_{p})  & 
C_{\ast}(X, B)\otimes H_{\ast}(F_{pq})  \\
};

\path[->, thick, font=\scriptsize]
(m-1-1)
edge node[anchor=east] {$\beta_{q}$} (m-2-1)
edge node[auto] {$q_{\ast}$} (m-3-2)
(m-1-3)
edge node[anchor=south] {$\simeq$} (m-1-1)
edge node[anchor=south east] {$q_{\ast}$} (m-3-2)
edge node[auto] {$\beta_{pq}$} (m-4-3)
(m-2-1)
edge node[anchor= east] {$\simeq$} (m-3-1)
(m-3-1)
edge node[anchor= east] {$\beta_{q}\otimes \id$} (m-4-1)
(m-3-2)
edge node[anchor=south] {$\cong$} (m-3-1)
edge node[auto] {$\beta_{q}$} (m-4-2)
(m-4-2)
edge node[anchor=south] {$\cong$} (m-4-1)
(m-4-3)
edge node[anchor=south] {$\ \ \ q^{\infty}$} (m-4-2)
(m-5-2)
edge node[auto] {$\cong$} (m-5-1)
(m-5-3)
edge node[auto] {$\id\otimes q_{\ast}$} (m-5-2)
;
\path[solid,  thick, font=\scriptsize]  
(m-4-1) 
edge node[anchor=east] {additivity} (m-5-1)
(m-4-3) 
edge node[anchor=west] {additivity} (m-5-3)
; 
\node[font=\Large] at (0.5,2.2) {\ding{192}};
\node[font=\Large] at (-1.4,0.9) {\ding{193}};
\node[font=\Large] at (-1.4, -0.8) {\ding{194}};
\node[font=\Large] at (2.5, 0.2) {\ding{195}};
\node[font=\Large] at (0.5,-3) {\ding{196}};
\end{tikzpicture}
\end{equation*}
Each vertex of this diagram represents a functor $\Rfd(B)\to \Chfd(\Q)$. The edges represent 
natural weak equivalences, with the exception of the lowest vertical edges where the passage 
between functors is obtained using additivity. The outer edges of this diagram correspond to
the homotopies in the diagram (\ref{PRODHOMOT2_MAPS}). Since 
$\tilde H_{\ast}(F_{q}) = 0$ the twisted tensor product  of the fibration $q^{\ast}p^{\ast}X \to p^{\ast}X$
is just  the untwisted tensor product $C_{\ast}(p^{\ast}X, E)\otimes H_{\ast}(F_{q})$. 
In effect the homotopy $\eta^{h}_{q}\circ (A(p^{!})\times \id_{I})$ in (\ref{PRODHOMOT2_MAPS})
is induced simply by the the quasi-isomorphisms 
$$
C_{\ast}(q^{\ast}p^{\ast}X, D)
\overset{\beta_{q}}{\lra} 
\bC(p^{\ast}X, E)\otimes H_{\ast}(F_{q})
\overset{\simeq}{\lra} C_{\ast}(p^{\ast}X, E)\otimes H_{\ast}(F_{q})$$ 
without using additivity. 

In order to show that the diagram (\ref{PRODHOMOT2_MAPS}) can be filled by a homotopy of 
homotopies it is enough to show that each of the subdiagrams in the above diagram of functors 
can be filled by a homotopy of homotopies. In the case of the subdiagrams (1)--(4) such homotopies 
of homotopies exist since subdiagrams (1) and (3) commute strictly and subdiagrams (2) and (4)
commute up to natural chain homotopies. Homotopy commutativity  of the subdiagram (2) follows 
from \cite[(7.4)]{Brown}.
The subdiagram (4) is homotopy commutative since it is obtained by applying 
Propositon \ref{MAPSFIBR-PROP} to the map of fibrations 
\begin{equation*}
\begin{tikzpicture}[baseline=(current bounding box.center)]
\matrix (m) 
[matrix of math nodes, row sep=2em, column sep=2em, text height=1.5ex, text depth=0.25ex]
{
D & &  E \\
& B &  \\
};
\path[->, thick, font=\scriptsize]
(m-1-1)
edge node[auto] {$q$} (m-1-3)
edge node[anchor= north east] {$pq$} (m-2-2)
(m-1-3)
edge node[anchor= north west] {$p$}  (m-2-2)
; 
\end{tikzpicture}
\end{equation*}
Finally, the subdiagram (5) can be filled by a homotopy of homotopies since the maps 
$$C_{\ast}(X, B)\otimes_{(pq)_{X}}H_{\ast}(F_{pq}) \overset{q^{\infty}}{\lra}
C_{\ast}(X, B)\otimes_{p_{X}}H_{\ast}(F_{p}) \overset{\cong}{\lra}
C_{\ast}(X, B)\otimes_{p_{X}}H_{\ast}(F_{p})\otimes H_{\ast}(F_{q})
$$
preserve the homological filtrations and induce isomorphisms on the filtration quotients.
\end{proof}

\begin{lemma}
\label{SECCOMPPROD-LEMMA}
Let  $p\colon E\to B$ be a unipotent fibration and let  $q\colon E\times S^{0}\to E$ be  
the product fibration with fiber $S^{0}$. Then 
$$\wh((pq)^{!})\simeq \wh(q^{!})\circ \wh(p^{!})$$
\end{lemma}

\begin{proof}
The basic outline of our argument is the same as in the proof of Lemma \ref{SECCOMPRATTRIV-LEMMA}. 
The same construction as in the proof of that lemma gives a homotopy  $H_{A}(p, q)$ 
between the maps $A((pq)^{!})$ and $A(p^{!})A(q^{!})$. Next, we need a homotopy $H_{K}$
between the maps $\chi(F_{p}\times S^{0})$ and $\chi(S^{0})\chi(F_{p})$. The first of these maps 
comes from the functor $C\mapsto C\otimes H_{\ast}(F_{p}\times S^{0})$ while the second map is induced 
by the functor $C\mapsto C\otimes H_{\ast}(F_{p})\otimes H_{\ast}(S^{0})$. The isomorphism 
$H(F_{p})\otimes H_{\ast}(S^{0}) \cong H_{\ast}(F_{p}\times S^{0})$
induces a natural isomorphism of the above functors, which in turn defines the homotopy 
$H_{K}$. 

As the result of these constructions we obtain a diagram of homotopies (\ref{PRODHOMOT2_MAPS})
(with $F_{q} = S^{0}$ and $F_{pq} = F_{p}\times S^{0}$). 
The remaining step is to show that this diagram can be filled by a homotopy of homotopies. 
Existence of such a homotopy of homotopies can be verified in a straightforward manner 
using the fact that for $X\in \Rfd(B)$ we have  a commutative diagram
\begin{equation*}
\begin{tikzpicture}[baseline=(current bounding box.center)]
\matrix (m) 
[matrix of math nodes, row sep=3em, column sep=2em, text height=1.5ex, text depth=0.25ex]
{
C_{\ast}(q^{\ast}p^{\ast}X, E\times S^{0}) 
&  C(p^{\ast}X, E)\otimes H_{\ast}(S^{0}) \\
\bC(X, B)\otimes_{(pq)_{X}} H(E\times S^{0}) 
& (C(X, B)\otimes_{p_{X}} H_{\ast}(F_{p}))\otimes H_{\ast}(S^{0})  \\
};
\path[->, thick, font=\scriptsize]
(m-1-1)
edge node[auto] {$\cong$} (m-1-2)
edge node[anchor=east] {$\beta_{pq}$} (m-2-1)
(m-1-2)
edge node[anchor=west] {$\beta_{p}\otimes\id$}  (m-2-2)
(m-2-1) 
edge  node[anchor=south] {$\cong$} (m-2-2)
; 
\end{tikzpicture}
\end{equation*}
\end{proof}

\begin{igusasarg}
\label{IGUSAARG NN}
In \cite{IgusaAx} Igusa used Lemma \ref{IGUSA-LEMMA}  to show that 
two higher torsion invariants of unipotent bundles coincide. We will adapt this 
argument to prove Theorems \ref{MAIN2 THM} and \ref{HMAIN2-THM}. The main 
idea of Igusa's  proof is to define a ``difference torsion'', which is new invariant that 
measures the difference between the given two invariants of bundles, and then to 
show that this difference torsion vanishes for all unipotent bundles. 
The essential properties of Igusa's difference torsion are encapsulated in 
Definition \ref{DIFFINV-DEF}. Proposition \ref{IGUSASARG-PROP} spells out the conditions 
that guarantee its vanishing. 
\end{igusasarg}

\begin{definition}
\label{DIFFINV-DEF}
Let $B$ be a  space of the homotopy type of a finite CW-complex, 
and let $\Lambda$ be an abelian group. 
An additive homotopy $B$-invariant of unipotent fibrations with values in $\Lambda$ is 
an assignment $\Phi$ 
that associates to each unipotent fibration  $p\colon E\to B$ an element $\Phi(p)\in \Lambda$
and  that satisfies the following conditions:
\begin{itemize}
\item[] \textbf{(additivity)}\ \ \  given maps of fibrations over $B$
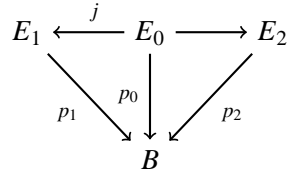
\begin{equation*}
\begin{tikzpicture}[baseline=(current bounding box.center)]
\matrix (m) 
[matrix of math nodes, row sep=3em, column sep=2.5em, text height=1.5ex, text depth=0.25ex]
{
E_{1} & E_{0}&  E_{2} \\
& B &  \\
};
\path[->, thick, font=\scriptsize]
(m-1-1)
edge node[anchor= north east] {$p_{1}$} (m-2-2)
(m-1-2)
edge node[anchor=  east] {$p_{0}$} (m-2-2)
edge node[anchor= south] {$\ \ j$} (m-1-1)
edge  (m-1-3)
(m-1-3)
edge node[anchor= north west] {$p_{2}$}  (m-2-2)
; 
\end{tikzpicture}
\end{equation*}
where $p_{i}$ is a unipotent fibration for $i=0, 1, 2$, and $j$ is a cofibration we have 
$$\Phi(p_{1}\cup_{p_{0}}p_{2}) = \Phi(p_{1})+ \Phi(p_{2})-\Phi(p_{0})$$
\item[] \textbf{(homotopy invariance)} if unipotent fibrations 
$p_{i}\colon E_{i}\to B$ ($i=1, 2$) are fiberwise homotopy equivalent then 
$\Phi(p_{1})= \Phi(p_{2})$.
\end{itemize}
\end{definition}

\begin{proposition}
\label{IGUSASARG-PROP}
Let $\Phi$ be an additive homotopy $B$-invariant with values in $\Lambda$. Assume that 
\begin{itemize}
\item $\Phi(p)=0$ for the product fibration $p\colon B\times S^{0} \to B$;
\item $\Phi(p)=0$ if the map $p\colon E\to B$ is a rational homotopy equivalence. 
\end{itemize}
Then $\Phi(p) = 0$ for all unipotent fibrations $p\colon E\to B$. 
\end{proposition}

\begin{proof}
First notice that $\Phi(p) = 0$ if $p\colon B\times F \to B$ is a product fibration where 
$F$ is either a disc or a sphere. Indeed, in the first case $p$ is a rational homotopy 
equivalence. If $F$ is a sphere we can argue inductively starting with $F=S^{0}$
and using additivity.   

Next,  let $p\colon E\to B$ be an arbitrary unipotent fibration. 
Applying Lemma \ref{IGUSA-LEMMA} to $p$ we obtain a sequence of fibrations
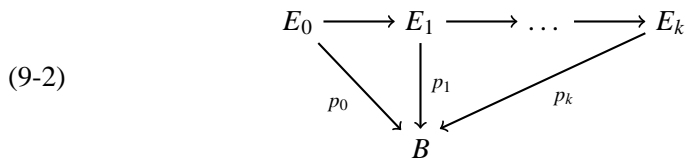
\begin{equation}
\label{IGUSAARG-EQ}
\begin{tikzpicture}[baseline=(current bounding box.center)]
\matrix (m) 
[matrix of math nodes, row sep=3em, column sep=2.5em, text height=1.5ex, text depth=0.25ex]
{
E_{0}& E_{1} & \dots & E_{k} \\
& B &   &  \\
};
\path[->, thick, font=\scriptsize]
(m-1-1)
edge  (m-1-2)
edge node[anchor=north east] {$p_{0}$} (m-2-2)
(m-1-2)
edge  (m-1-3)
edge node[anchor= west] {$p_{1}$} (m-2-2)
(m-1-3)
edge  (m-1-4)
(m-1-4)
edge node[anchor=north west] {$p_{k}$} (m-2-2.north east)
; 
\end{tikzpicture}
\end{equation}
Since $p_{k}$ is a rational homotopy equivalence we have $\Phi(p_{k})=0$.
Next, since  $p_{i}$ is obtained as a pushout  of $p_{i-1}$ and product fibrations with disc 
and sphere fibers we can use additivity of $\Phi$ to  get that 
$\Phi(p_{i}) = \Phi(p_{i-1})$. This gives 
$$0 = \Phi(p_{k}) = \Phi(p_{k-1}) = \dots = \Phi(p_{0})$$  
Finally, since $p_{0}$ is an $n$-fold fiberwise suspension of $p$ we can use additivity of 
$\Phi$ again to get 
$$\Phi(p) = (-1)^{n}\Phi(p_{0}) = 0$$
\end{proof}

\begin{note}
\label{ADDBUNDLEINV NOTE}
Igusa uses a  variant of Proposition \ref{IGUSASARG-PROP} that will be also 
useful to us later on. Namely, for a smooth compact manifold  $B$ consider an assignment 
$\Phi$ that satisfies additivity and homotopy invariance properties as in 
Definition \ref{DIFFINV-DEF}, but is defined only for unipotent bundles over $B$. Then 
the statement of  Proposition \ref{IGUSASARG-PROP} still holds: if $\Phi$ vanishes on 
the product bundle $B\times S^{0}\to B $ and on bundles that are given by a rational homotopy 
equivalence then it vanishes on all unipotent bundles. The proof of this fact is essentially 
the same as the proof of Proposition \ref{IGUSASARG-PROP}  with two additional observations:

\textbullet\  the homotopy invariance of $\Phi$ lets us define this invariant on all 
smoothable unipotent fibrations over  $B$, i.e. all  fibrations that are 
fiberwise homotopy equivalent to a unipotent bundle. 

\textbullet\ if $p\colon E\to B$ is a unipotent bundle then all fibrations appearing in 
the diagram (\ref{IGUSAARG-EQ}) are smoothable \cite[Lemma 8.6]{IgusaAx}. 

\end{note}

We are now ready to give 

\begin{proof}[Proof of Theorem  \ref{HMAIN2-THM}]
Let  $p\colon E\to B$,  $q\colon D\to E$ be unipotent fibrations. We have 
commutative diagrams 
\begin{equation*}
\begin{tikzpicture}[baseline=(current bounding box.center)]
\matrix (m) 
[matrix of math nodes, row sep=3em, column sep=5em, text height=1.5ex, text depth=0.5ex]
{
\wh(B) & \wh(D) \\
A(B) &  A(D) \\
};
\path[->, thick, font=\scriptsize]
(m-1-1)
edge node[auto] {$\wh((pq)^{!})$} (m-1-2)
edge node[anchor=east] {$i_{B}$} (m-2-1)
(m-1-2)
edge node[anchor=west] {$i_{D}$}  (m-2-2)
(m-2-1) 
edge  node[anchor=north] {$A((pq)^{!})$} (m-2-2)
; 
\matrix (m) 
[xshift=70mm, matrix of math nodes, row sep=3em, column sep=5em, text height=1.5ex, text depth=0.5ex]
{
\wh(B) & \wh(D) \\
A(B) &  A(D) \\
};
\path[->, thick, font=\scriptsize]
(m-1-1)
edge node[auto] {$\wh(q^{!})\wh(p^{!})$} (m-1-2)
edge node[anchor=east] {$i_{B}$} (m-2-1)
(m-1-2)
edge node[anchor=west] {$i_{D}$}  (m-2-2)
(m-2-1) 
edge  node[anchor=north] {$A(q^{!})A(p^{!})$} (m-2-2)
; 
\end{tikzpicture}
\end{equation*}
In the same way as in the proof of Lemma \ref{SECCOMPRATTRIV-LEMMA}
we  construct a homotopy $H_{A}(p, q)$ between the maps $A((pq)^{!})$ and $A(q^{!})A(p^{!})$. 
By  (\ref{HOMOTOBSTR NN}) this data defines a map $\varphi(p, q)\colon \wh(B)\to \Omega K(\Q)$
such that $[\varphi(p, q)]=0$  if and only if $H_{A}$ admits a lift to a homotopy 
between $\wh((pq)^{!})$ and $\wh(q^{!})\wh(p^{!})$. 

Fix a fibration $p\colon E\to B$. Let $\Phi_{p}$ be the assignment that associates to 
a unipotent fibration $q\colon D\to E$ the homotopy class $[\varphi(p, q)]$. 
We claim that $\Phi_{p}$ is an additive  homotopy $E$-invariant with values in the group 
$[\wh(B), \Omega K(\Q)]$ (\ref{DIFFINV-DEF}).  In order to  verify homotopy invariance of $\Phi_{p}$
assume that we have a fiberwise homotopy equivalence
\begin{equation*}
\begin{tikzpicture}[baseline=(current bounding box.center)]
\matrix (m) 
[matrix of math nodes, row sep=2em, column sep=2em, text height=1.5ex, text depth=0.25ex]
{
D_{1} & &  D_{2} \\
& E &  \\
};
\path[->, thick, font=\scriptsize]
(m-1-1)
edge node[auto] {$f$} (m-1-3)
edge node[anchor= north east] {$q_{1}$} (m-2-2)
(m-1-3)
edge node[anchor= north west] {$q_{2}$}  (m-2-2)
; 
\end{tikzpicture}
\end{equation*}
We need to check that the maps $\varphi(p, q_{1})$ and $\varphi(p, q_{2})$ are homotopic. 
By  Proposition \ref{HOMOT INV PROP} we can construct  a homotopy 
$$H^{\rm{Wh}(E)}_{f}\colon \wh(E)\times I\to \wh(D_{2})$$ 
between the maps 
$f_{\ast}\wh(q_{1}^{!})$ and $\wh(q_{2}^{!})$.  As a consequence the map 
$$\bar{H}^{\rm{Wh}(E)}_{f} := H^{\rm{Wh}(E)}_{f}\circ (\wh(p^{!})\times \id_{I})$$
is a homotopy between 
$f_{\ast}\wh(q_{1}^{!})\wh(p^{!})$ and $\wh(q_{2}^{!})\wh(p^{!})$. On the other hand 
$f$ is also a fiberwise homotopy equivalence of  the fibrations $pq_{1}$ and $pq_{2}$, 
so we have a homotopy 
$$H^{\rm{Wh}(B)}_{f}\colon \wh(B)\times I\to \wh(D_{2})$$ 
between the maps $f_{\ast}\wh((pq_{1})^{!})$ and $\wh((pq_{2})^{!})$.
By the proof of Proposition \ref{HOMOT INV PROP} we also have a homotopy
$$H^{A(B)}_{f}\colon A(B)\times I \to A(D_{2})$$ 
between  the maps $f_{\ast}A((pq_{1})^{!})$ and $A((pq_{2})^{!})$
as well as a homotopy 
$$\bar{H}^{A(E)}_{f} := H^{A(E)}_{f}\circ (A(p^{!})\times \id_{I})$$
between $f_{\ast}A(q^{!}_{1})A(p^{!})$ and $A(q^{!}_{2})A(p^{!})$. 
All these homotopies fit into commutative diagrams 
\begin{equation*}
\begin{tikzpicture}[baseline=(current bounding box.center)]
\matrix (m) 
[matrix of math nodes, row sep=3em, column sep=4em, text height=1.5ex, text depth=0.5ex]
{
\wh(B) \times I& \wh(D_{2}) \\
A(B) \times I &  A(D_{2}) \\
};
\path[->, thick, font=\scriptsize]
(m-1-1)
edge node[auto] {$\ \ \ \ \ \ H^{\rm{Wh}(B)}_{f}$} (m-1-2)
edge node[anchor=east] {$i_{B}\times \id_{I}$} (m-2-1)
(m-1-2)
edge node[anchor=west] {$i_{D_{2}}$}  (m-2-2)
(m-2-1) 
edge  node[anchor=north] {$H^{A(B)}_{f}$} (m-2-2)
; 
\matrix (m) 
[xshift=70mm, matrix of math nodes, row sep=3em, column sep=4em, text height=1.5ex, text depth=0.5ex]
{
\wh(B)\times I & \wh(D) \\
A(B)\times I &  A(D) \\
};
\path[->, thick, font=\scriptsize]
(m-1-1)
edge node[auto] {$\ \ \ \ \ \ \bar{H}^{\rm{Wh}(E)}_{f}$} (m-1-2)
edge node[anchor=east] {$i_{B}\times \id_{I}$} (m-2-1)
(m-1-2)
edge node[anchor=west] {$i_{D_{2}}$}  (m-2-2)
(m-2-1) 
edge  node[anchor=north] {$\bar{H}^{A(E)}_{f}$} (m-2-2)
; 
\end{tikzpicture}
\end{equation*}
Consider the diagram 
\begin{equation}
\label{HOMINVDIFF EQ}
\begin{tikzpicture}[baseline=(current bounding box.center)]
\matrix (m) 
[matrix of math nodes, row sep=3em, column sep=5em, text height=1.5ex, text depth=0.5ex]
{
f_{\ast} A((pq_{1})^{!}) & f_{\ast} A(q_{1}^{!})A(p^{!}) \\
A((pq_{2})^{!}) &  A(q_{2}^{!})A(p^{!})  \\
};
\path[solid, thick, font=\scriptsize]
(m-1-1)
edge node[auto] {$f_{\ast}H_{A}(p,q_{1})$} (m-1-2)
edge node[anchor=east] {$H^{A(B)}_{f}$} (m-2-1)
(m-1-2)
edge node[anchor=west] {$\bar{H}^{A(E)}_{f}$}  (m-2-2)
(m-2-1) 
edge  node[anchor=north] {$H_{A}(p,q_{2})$} (m-2-2)
; 
\end{tikzpicture}
\end{equation}
Each vertex of this diagram represents a map $A(B)\to A(D_{2})$ and edges 
represent homotopies of such maps. This diagram is induced by  a diagram of functors 
$\Rfd(B)\to \Rfd(D_{2})$ and natural equivalences of such functors. It is straightforward 
to check that this underlying diagram of functors commutes. This implies that the diagram 
(\ref{HOMINVDIFF EQ}) can be filled by a homotopy of homotopies. 
This homotopy of homotopies can be interpreted as a homotopy between $H^{A(B)}_{f}$ and 
$\bar{H}^{A(E)}_{f}$. By  (\ref{HOMOTOBSTR NN}) this homotopy defines a map 
$\wh(B)\times I \to \Omega K(\Q)$. One can check that this map determines a homotopy between 
$\varphi(p, q_{1})$ and $\varphi(p, q_{2})$. 
Additivity of $\Phi_{p}$ can be verified in a similar way, using additivity of secondary transfers.

By Lemma \ref{SECCOMPRATTRIV-LEMMA} and Lemma \ref{SECCOMPPROD-LEMMA}
we have $\Phi_{p}(q) = 0$ if $q$ is a product fibration or a rational homotopy equivalence. 
Proposition  \ref{IGUSASARG-PROP} implies then that $[\varphi(p, q)] = \Phi_{p}(q)=0$ 
for any unipotent fibration $q\colon D\to E$.
 \end{proof}

\begin{proof}[Proof of Theorem \ref{MAIN2 THM}] 
Let $p\colon E\to B$ and $q\colon D\to E$ be unipotent bundles. 
Recall (\ref{TRANSFER CONSTR NN}) that  the map  $Q(p!)$ is induced by functors 
between categories of partitions that assign to a partition $P\subseteq B\times I$
the partition $(p\times \id)^{-1}(P)\subseteq E\times I$. Since 
$ (pq\times \id)^{-1}(P) = (q\times \id)^{-1}(p\times \id)^{-1}(P)$ we obtain a homotopy $H_{Q}(p, q)$ between 
between the maps  $Q((pq)^{!})$ and $Q(q^{!})Q(p^{!})$. By (\ref{HOMOTOBSTR NN}) the maps 
$\wh[s](q^{!})\wh[s](p^{!})$, $\wh[s]((pq)^{!})$ and the homotopy $H_{Q}(p, q)$ define a map 
$\psi(p, q)\colon \wh[s](B)\to \Omega K(\Q)$ such that $\wh[s](q^{!})\wh[s](p^{!})\simeq \wh((pq)^{!})$ if $[\psi(p, q)] = 0$. 

It remains to show that $[\psi(p, q)]=0$ for all unipotent bundles $p$, $q$. 
Recall (\ref{MAIN1HOMOTRED-NN}) that by $\mu_{p}$ we denoted the homotopy between the maps
$A(p^{!})a_{B}$ and $a_{E}Q(p^{!})$. Consider the diagram  
\begin{equation}
\label{HAHQ-EQ}
\begin{tikzpicture}[baseline=(current bounding box.center)]
\matrix (m) 
[matrix of math nodes, row sep=3em, column sep=6em, text height=1.5ex, text depth=0.25ex]
{
a_{D}Q((pq)^{!}) &  a_{D}Q(q^{!})Q(p^{!}) \\
 & A(q^{!})a_{E}Q(p^{!})\\
A((pq)^{!})a_{B} &  A(q^{!})A(p^{!})a_{B} \\
};
\path[->, thick, font=\scriptsize]
(m-1-1)
edge node[auto] {$a_{D}H_{Q}(p,q)$} (m-1-2)
edge node[anchor=east] {$\mu_{pq}$} (m-3-1)
(m-1-2)
edge node[anchor=west] {$\mu_{q}\circ (Q(p^{!})\times \id_{I})$}  (m-2-2)
(m-2-2) 
edge  node[anchor=west] {$H_{A}(p, q)\circ (a_{B}\times \id_{I})$} (m-3-2)
(m-3-1) 
edge  node[anchor=north] {$H_{A}(p, q)\circ (a_{B}\times \id_{I})$} (m-3-2)
; 
\end{tikzpicture}
\end{equation}
Vertices of this diagram corresponds to a map $Q(B_{+}) \to A(D)$ and edges 
correspond to homotopies of such maps. Each of these homotopies is induced by 
natural isomorphisms between retractive spaces over $D$ obtained from partitions $P\subseteq B\times I$:
\begin{equation*}
\begin{tikzpicture}[baseline=(current bounding box.center)]
\matrix (m) 
[matrix of math nodes, row sep=2em, column sep=3em, text height=1.5ex, text depth=0.25ex]
{
(pq\times \id)^{-1}(P) &  (q\times \id)^{-1}(p\times\id)^{-1}(P) \\
 & q^{\ast}(p\times \id)^{-1}(P)\\
(pq)^{\ast}P &  q^{\ast}p^{\ast}P \\
};
\path[->, thick, font=\scriptsize]
(m-1-1)
edge  (m-1-2)
edge (m-3-1)
(m-1-2)
edge   (m-2-2)
(m-2-2) 
edge   (m-3-2)
(m-3-1) 
edge  (m-3-2)
; 
\end{tikzpicture}
\end{equation*} 
Since this diagram commutes we obtain a homotopy of homotopies filling 
the diagram (\ref{HAHQ-EQ}).

Let $b_{B}\colon \wh[s](B)\to \wh(B)$ be the map induced by the map of fibrations
\begin{equation*}
\begin{tikzpicture}[baseline=(current bounding box.center)]
\matrix (m) 
[matrix of math nodes, row sep=2em, column sep=1em, text height=1.5ex, text depth=0.5ex]
{
\wh[s](B)&  & \wh(B) \\
Q(B_{+}) & &  A(B) \\
& K(\Q) & \\
};
\path[->, thick, font=\scriptsize]
(m-1-1)
edge node[auto] {$b_{B}$} (m-1-3)
edge  (m-2-1)
(m-1-3)
edge  (m-2-3)
(m-2-1)
edge node[auto] {$a_{B}$} (m-2-3)
edge node[anchor=north east] {$\lambda_{B}$} (m-3-2)
(m-2-3)
edge node[auto] {$\lambda^{h}_{B}$} (m-3-2)

; 
\end{tikzpicture}
\end{equation*}
The construction of the map $\psi(p, q)$ described in (\ref{HOMOTOBSTR NN})
combined with existence of a homotopy of homotopies in the diagram (\ref{HAHQ-EQ}) gives
$$\psi(p,q) \simeq \varphi(p, q) b_{B} $$
where $\varphi(p, q)$ is the map defined in the proof of Theorem \ref{HMAIN2-THM}. 
In that proof we showed that $[\varphi(p,q)] = 0$, so also $[\psi(p, q)] = 0$.
\end{proof}

%%%%%%%%%%%%%%%%%%%%%%%%%%%%
% SECONDARY TRANSFER AND SMOOTH TORSION
%%%%%%%%%%%%%%%%%%%%%%%%%%%%

\section{Secondary transfer and smooth torsion}
\label{TORSION SEC}

In \cite{BDW} and \cite{BDKW} (joint with B. Williams and J. Klein) we
described a homotopy theoretical construction of the smooth torsion of unipotent bundles 
and showed that it defines characteristic classes which coincide 
with the higher torsion invariants of Igusa and Klein. The construction of the smooth torsion 
of a bundle $p\colon B\to E$ proceeds as follows. Let $\eta_{B}\colon B\to Q(B_{+})$ 
denote the coaugmentation map. By  \cite[Theorem 6.7]{BDW} the map 
$$\lambda_{E}Q(p^{!})\eta_{B} \colon B \to K(\Q)$$
is homotopic via a preferred homotopy $\omega_{p}$
to a constant map. This defines  a map $\tau^{s}(p)\colon B \to \wh[s](E)$
which is a lift of $Q(p^{!})\eta_{B}$. The map $\tau^{s}(p)$ is the smooth torsion of the bundle $p$. 

The secondary transfer of unipotent bundles described in this paper  can be used
to construct another map  $\bar{\tau}^{s}(p)\colon B\to \wh[s](E)$. Namely, since 
the identity map $\id_{B}\colon B\to B$ can be considered as a unipotent bundle, 
it defines the smooth torsion $\tau^{s}(\id_{B})\colon B\to \wh[s](B)$. We set
$$\bar{\tau}^{s}(p) := \wh[s](p^{!})\tau^{s}(\id_{B})$$

Our final goal in this paper is to prove  Theorem \ref{MAIN3 THM} which says 
that for any unipotent bundle $p$ the maps $\tau^{s}(p)$ and $\bar{\tau}^{s}(p)$ are homotopic.  
We will also show that as a consequence the statement of Theorem  \ref{THM_TRAX} holds:
for any pair of composable unipotent bundles $p$ and $q$ the higher torsion  cohomology 
classes of $p$, $q$ and $pq$ are related by  the formula (\ref{EQ-TRAX}).

The proof of Theorem \ref{MAIN3 THM} will use the same scheme as the proof of the composition 
formula of unipotent fibrations (Proposition \ref{HMAIN2-THM}). We will first show that 
this theorem holds when $p$ is either a product bundle or it is a rational homotopy equivalence, 
and then we will use Igusa's argument (\ref{IGUSAARG NN})
to extend this result to all unipotent bundles.

\begin{lemma}
\label{TORSIONRATPROD LEMMA}
Let $p\colon E\to B$ be a unipotent bundle. If $p$ is either the product bundle 
with fiber $S^{0}$ or a rational homotopy equivalence then $\tau^{s}(p)\simeq \bar{\tau}^{s}(p)$. 
\end{lemma}

\begin{proof}
This follows essentially from \cite[\S 6]{BDKW}.  We proved there 
that if $p\colon E\to B$ is a unipotent bundle that satisfies the assumptions 
of the Leray-Hirsch theorem then  the quasi-isomorphisms 
$$C_{\ast}(p^{\ast}X) \overset{\simeq}{\lra} C_{\ast}(X)\otimes H_{\ast}(F_{p})$$
given for $X\in \Rfd(B)$ by that theorem define a map 
$$\wh[LH](p^{!})\colon \wh[s](B)\to \wh[s](E)$$
and that $\tau^{s}(p)\simeq \wh[LH](p^{!})\tau^{s}(\id_{B})$. It is straightforward to check that 
if $p$ is either a product bundle $B\times S^{0}\to B$ or a rational homotopy equivalence then 
$\wh[LH](p^{!})\simeq \wh[s](p^{!})$.  
\end{proof}

\begin{proof}[Proof of Theorem  \ref{MAIN3 THM}]
Both   $\tau^{s}(p)$ and $\bar{\tau}^{s}(p)$ are defined as lifts of the map 
$Q(p^{!})\eta_{B}$.
As a consequence they define a map 
$\varrho(p)\colon B \to \Omega K(\Q)$ such that the homotopy class of the composition
$$B \overset{\varrho(p)}{\lra} \Omega K(\Q) \lra \wh[s](E)$$
coincides with the element $[\tau^{s}(p)] - [\bar{\tau}^{s}(p)] \in [B, \wh[s](E)]$. It  suffices 
to show that $[\varrho(p)]= 0$ in $[B, \Omega K(\Q)]$. 

We claim that the assignment $p \mapsto [\varrho(p)]$ is an additive homotopy $B$-invariant 
of unipotent bundles with values in $[B, \Omega K(\Q)]$ (\ref{ADDBUNDLEINV NOTE}). 
Indeed, additivity of this assignment follows essentially from the additivity of the secondary 
transfer
(\ref{SADD THM}) and the additivity  of the smooth torsion
\cite[Theorem 5,1]{BDKW}.  Homotopy invariance can be verified
using the fact that  the construction of $\varrho(p)$ involves only  
chain complexes associated to  $p$.  
Using Lemma \ref{TORSIONRATPROD LEMMA} we obtain that $[\varrho(p)]=0$ 
if $p$ is either a product bundle or a rational homotopy equivalence. 
By Proposition \ref{IGUSASARG-PROP} and (\ref{ADDBUNDLEINV NOTE}) 
we get  then that $[\varrho(p)] = 0$ for all unipotent bundles $p$.

\end{proof}

\begin{proof}[Proof of Theorem \ref{THM_TRAX}]
Combining Theorems  \ref{MAIN3 THM} and  \ref{MAIN2 THM} we obtain the statement of 
Corollary \ref{MAIN3 COR}: for any unipotent bundles $p\colon E\to B$ and $q\colon D\to E$ 
there is a homotopy
$$\tau^{s}(pq) \simeq \wh[s](q^{!})\tau^{s}(p)$$
In \cite[Theorem 7.1]{BDKW} an analogous decomposition of  the smooth torsion of $pq$
(in the case where $q$ is a Leray-Hirsch bundle) was the main ingredient in the proof 
of the fact that the cohomological torsion of $p$ and $q$ satisfies the  formula  (\ref{EQ-TRAX}). 
The same argument can be now used to obtain the formula  (\ref{EQ-TRAX}) for arbitrary 
unipotent bundles $p$ and $q$. 
\end{proof}

%%%%%%%%%%%%%%%%%%%%%%%%%%%%
% APPENDIX: MAPS OF HOMOTOPY FIBERS
%%%%%%%%%%%%%%%%%%%%%%%%%%%%

\section{Appendix A: Maps of homotopy fibers}
\label{HOMOTFIB APP}

Multiple arguments in this paper involve constructions of maps between homotopy fibers 
as well as constructions of homotopies between such maps. We summarize here 
the basic scheme of such constructions. 

\begin{homotfibmap} 
\label{HOMOTFIBMAP NN}
For a space $X$ with a basepoint $x_{0}$ let $P_{\!\!\! x_{0}}X$ denote
the space of paths in $X$ that start at $x_{0}$. By the homotopy fiber of a map  
$p\colon Y \to X$ over  $x_{0}$ we understand the standard construction
$$\hofib(p)_{x_{0}} := \{(\omega, y) \in P_{\!\!\! x_{0}}X\times Y \ | \ \omega(1) = p(y) \}$$ 
We will denote by  $i_{Y}\colon \hofib(p)_{x_{0}}\to Y$ the map  given by $i_{Y}(\omega, y) = y$. 

 Assume  that we have a diagram 
\begin{equation*}
\begin{tikzpicture}
\matrix (m) 
[matrix of math nodes, row sep=2.5em, column sep=2.5em, text height=1.5ex, text depth=0.25ex]
{
Y & Y' \\
X & X' \\ 
};
\path[->, thick, font=\scriptsize]
(m-1-1)
edge node[auto] {$\tilde f$} (m-1-2)
edge node[anchor=east] {$p$} (m-2-1)
(m-1-2) 
edge node[auto] {$p'$} (m-2-2)
(m-2-1)
edge node[auto] {$f$} (m-2-2)
; 
\end{tikzpicture}
\end{equation*}
such that $f(x_{0}) = x'_{0}$. Given a homotopy $h$ from $fp$ to  $p'\tilde f$ we obtain 
a map $\tilde{\tilde{f}}\colon \hofib(p)_{x_{0}}\to \hofib(p')_{x'_{0}}$ given by 
$$\tilde{\tilde{f}}(\omega, y) := (f\omega \ast h_{y}, \ \tilde{f}(y))$$ 
Here $h_{y}$ denotes the path in $X'$ defined by $h_{y}(t) = h(y, t)$, and $\ast$
indicates concatenation of paths. We have 
$$i_{Y'}\tilde{\tilde f} = \tilde{f} i_{Y}$$
\end{homotfibmap}

\begin{homotfibhomot}
\label{HOMOTFIBHOMOT NN}
Assume   that we have two diagrams:
\begin{equation}
\label{2HOMOTCOMMSQS EQ}
\begin{tikzpicture}[baseline=(current bounding box.center)]
\matrix (m) 
[matrix of math nodes, row sep=2.5em, column sep=2.5em, text height=1.5ex, text depth=0.25ex]
{
Y & Y' \\
X & X' \\ 
};
\path[->, thick, font=\scriptsize]
(m-1-1)
edge node[auto] {$\tilde{f}_{0}$} (m-1-2)
edge node[anchor=east] {$p$} (m-2-1)
(m-1-2) 
edge node[auto] {$p'$} (m-2-2)
(m-2-1)
edge node[auto] {$f_{0}$} (m-2-2)
; 
\matrix (n) 
[xshift=40mm, matrix of math nodes, row sep=2.5em, column sep=2.5em, text height=1.5ex, text depth=0.25ex]
{
Y & Y' \\
X & X' \\ 
};
\path[->, thick, font=\scriptsize]
(n-1-1)
edge node[auto] {$\tilde{f}_{1}$} (n-1-2)
edge node[anchor=east] {$p$} (n-2-1)
(n-1-2) 
edge node[auto] {$p'$} (n-2-2)
(n-2-1)
edge node[auto] {$f_{1}$} (n-2-2)
; 
\end{tikzpicture}
\end{equation}
that commute up to homotopies $h_{0}$ and $h_{1}$, respectively, and such that 
$f_{0}(x_{0}) = f_{1}(x_{0}) = x'_{0}$. This gives  two 
maps of the homotopy fibers:
$$\tilde{\tilde{f}}_{0}, \tilde{\tilde{f}}_{1}\colon \hofib(p)_{x_{0}}\to \hofib(p')_{x'_{0}}$$
In order to obtain a homotopy between these maps it suffices to construct the following data:
\begin{list}{\labelitemi}{\leftmargin=0.5em}
\item[1)] a homotopy $\tilde{H} \colon Y\times I \to Y'$ between $\tilde{f}_{0}$ and $\tilde{f}_{1}$;
\item[2)] a basepoint preserving homotopy $H\colon X\times I \to X'$ between $f_{0}$ and $f_{1}$;
\item[3)] a homotopy of homotopies between the two homotopies $p'\tilde f_{0}\simeq f_{1}p$:
the one given  as a concatenation of  $p'\tilde{H}$ with $h_{1}$, and the one 
obtained by  concatenating $h_{0}$ with $H(p\times \id_{I})$:
\begin{equation*}
\begin{tikzpicture}[baseline=(current bounding box.center)]
\matrix (m) 
[matrix of math nodes, row sep=2.5em, column sep=5em, text height=1.5ex, text depth=0.25ex]
{
f_{0}p & f_{1}p \\ 
p'\tilde{f}_{0} &  p'\tilde{f}_{1}  \\
};
\path[->, thick, font=\scriptsize]
(m-2-1)
edge node[anchor= north] {$p'\tilde H$} (m-2-2)
edge node[anchor=east] {$h_{0}$} (m-1-1)
(m-2-2) 
edge node[anchor = west] {$h_{1}$} (m-1-2)
(m-1-1)
edge node[anchor=south] {$H(p\times \id_{I})$} (m-1-2)
; 
\end{tikzpicture}
\end{equation*}
\end{list}
Giving such a homotopy of homotopies is equivalent to giving a map
$$\Theta \colon Y \times I\times I \to X'$$
such that $\Theta |_{Y\times \{i\}\times I} = h_{i}$ for $i=0, 1$, 
$\Theta_{Y\times I \times \{0\}} = p'\tilde H$, 
and $\Theta_{Y\times I \times \{1\}} = H(p\times \id_{I})$.
The homotopy $\tilde{\tilde H}$ between $\tilde{\tilde{f}}_{1}$ and 
$\tilde{\tilde{f}}_{2}$ is defined by:
$$\tilde{\tilde{H}}((\omega, y), t) := ( H_{t}\omega \ast \Theta_{y,t} , \ \tilde{H}(y, t))$$ 
where $H_{t}\colon X\to X'$ is given by $H_{t}(x) = H(t, x)$ and $\Theta_{y, t}$ is the path 
in $X'$ given by $\Theta_{y, t}(s) = \Theta(y, t, s)$.
\end{homotfibhomot}

\begin{homotobstr}
\label{HOMOTOBSTR NN}
Assume again that we have two homotopy commutative squares 
(\ref{2HOMOTCOMMSQS EQ}) and that $\tilde{\tilde{f}}_{0}, \tilde{\tilde{f}}_{1}$
are the maps of homotopy fibers defined by these squares. Assume also that 
we have a homotopy $\tilde H$ between $\tilde{f}_{0}$ and $\tilde{f}_{1}$, and that we 
want  to determine whether there exists a homotopy $\tilde{\tilde H}$ between  
$\tilde{\tilde{f}}_{1}$ and $\tilde{\tilde{f}}_{2}$ that fits into a commutative  
diagram
\begin{equation}
\label{HTILTIL EQ}
\begin{tikzpicture}[baseline=(current bounding box.center)]
\matrix (m) 
[matrix of math nodes, row sep=2.5em, column sep=4em, text height=1.5ex, text depth=0.25ex]
{
\hofib(p)_{x_{0}}\times I  &\hofib(p')_{x'_{0}}  \\
Y\times I  & Y' \\ 
};
\path[->, thick, densely dashed, font=\scriptsize]
(m-1-1)
edge node[auto] {$\tilde{\tilde H}$} (m-1-2)
;
\path[->, thick, font=\scriptsize]
(m-1-1)
edge node[anchor=east] {$i_{Y}\times \id_{I}$} (m-2-1)
(m-1-2) 
edge node[auto] {$i_{Y'}$} (m-2-2)
(m-2-1)
edge node[anchor=north] {$\tilde H$} (m-2-2)
; 
\end{tikzpicture}
\end{equation}
In (\ref{HOMOTFIBHOMOT NN}) 
we described data that suffices to construct such ${\tilde{\tilde H}}$, but here we are interested 
in  a condition that is equivalent to the existence of this homotopy.  We can describe such a 
condition as follows.
\!\!\!\footnote{
See also \cite[Lemma 4.1]{PetersonEmery}
} 

Let $C_{p}\colon \hofib(p)_{x_{0}}\times I \to X$ be the map given by $C_{p}((\omega, y), t) = \omega(t)$. 
This is a homotopy between the constant map into $x_{0}$ and the map $pi_{Y}$. 
Consider the following diagram:
\begin{equation*}
\begin{tikzpicture}[baseline=(current bounding box.center)]
\matrix (m) 
[matrix of math nodes, row sep=2em, column sep=2em, text height=1.5ex, text depth=0.5ex]
{
p'\tilde{f}_{0}i_{Y} &  & p'\tilde{f}_{1}i_{Y} \\
f_{0}pi_{Y} & &  f_{1}pi_{Y'} \\
& x'_{0} & \\
};
\path[solid, thick, font=\scriptsize]
(m-1-1)
edge node[auto] {$p'\tilde{H}(i_{Y\times \id_{I}})$} (m-1-3)
edge node[anchor=east] {$h_{0}(i_{Y}\times \id_{I})$} (m-2-1)
(m-1-3) 
edge node[anchor=west] {$h_{1}(i_{Y'}\times \id_{I})$} (m-2-3)
(m-2-1)
edge node[anchor=north east] {$f_{0}C_{p}$} (m-3-2)
(m-2-3)
edge node[auto] {$f_{1}C_{p}$} (m-3-2)
; 
\end{tikzpicture}
\end{equation*}
Each vertex of this diagram represent a map $ \hofib(p)_{x_{0}} \to X'$ and  edges 
represent homotopies of such maps. Concatenating all homotopies appearing here 
we obtain a homotopy from the constant map into $x'_{0}$ to itself, or equivalently a map 
$$\varphi\colon \hofib(p)_{x_{0}} \to \Omega X'$$
It is straightforward to verify that the map $\varphi$ is contractible if and only if there exists a homotopy
$\tilde{\tilde H}$ such that the diagram (\ref{HTILTIL EQ}) commutes. In  other words 
the homotopy class  of $\varphi$ is an obstruction to lifting the homotopy $\tilde H$ 
to a homotopy defined on the level of the homotopy fibers. 
\end{homotobstr}

\textbf{Note.} Let  $p'\colon (Y', y'_{0})\to (X', x'_{0})$ be a map of infinite loop spaces 
where $x'_{0}$, $y'_{0}$ are the trivial elements in $X'$ and $Y'$. In this case the map
$\varphi$ has a simpler interpretation. Namely, let 
$j_{Y'}\colon \Omega X' \to \hofib(p')_{x'_{0}}$ be the map given by $j_{Y'}(\omega) = (\omega, y_{0})$. 
The set of homotopy classes $[\hofib(p)_{x_{0}}, \hofib(p')_{x'_{0}}]$ has a structure 
of an abelian group and we have:
$$j_{Y'\ast}[\varphi] = [\tilde{f}_{1}] - [\tilde{f}_{2}]$$

%%%%%%%%%%%%%%%%%%%%%%%%%%%%
%  BIBLIOGRAPHY
%%%%%%%%%%%%%%%%%%%%%%%%%%%%

\bibliographystyle{plain}
\bibliography{secondary}

\end{document}